\documentclass[leqno,12pt]{article}
\usepackage[T1]{fontenc}
\usepackage{amsmath}
\usepackage{amsthm}
\usepackage{amsxtra}
\usepackage{amsfonts}
\usepackage{amssymb}
\usepackage[margin=1.25in]{geometry}
\usepackage{color}
\usepackage{epsfig}
\usepackage{graphicx}
\usepackage{subfig}
\usepackage{cite}
\usepackage{color}
\usepackage{array}
\usepackage{booktabs}
\usepackage{fullpage}

\renewcommand{\div}{\text{div}}

\newcommand{\R}{\mathbb{R}}

\newcommand{\vb}[1]{\mathbf{#1}}

\newcommand{\eps}{\varepsilon}

\renewcommand{\phi}{\varphi}

\def\XXint#1#2#3{{\setbox0=\hbox{$#1{#2#3}{\int}$ }
\vcenter{\hbox{$#2#3$ }}\kern-.6\wd0}}

\newtheorem{theorem}{Theorem}
\newtheorem*{theorem*}{Theorem}
\theoremstyle{plain}

\newtheorem{lemma}{Lemma}

\title{
Accelerated PDE's for efficient solution of regularized inversion problems
}
\author{Minas Benyamin\thanks{School of Electrical and Computer Engineering, Georgia Institute of Technology. ({\tt minasbenyamin@gatech.edu})} \and Jeff Calder\thanks{School of Mathematics, University of Minnesota. ({\tt jcalder@umn.edu})} \and Ganesh Sundaramoorthi\thanks{United Technologies Research Center, East Hartford, CT ({\tt ganesh.sun@gmail.com})} \and Anthony Yezzi\thanks{School of Electrical and Computer Engineering, Georgia Institute of Technology.  ({\tt anthony.yezzi@ece.gatech.edu})}}

\begin{document}
\maketitle
\begin{abstract}
We further develop a new framework, called \emph{PDE Acceleration},
by applying it to calculus of variations problems defined for general
functions on $\R^n$,
obtaining efficient numerical algorithms to solve the resulting 
class of optimization problems based on simple discretizations of their 
corresponding accelerated PDE's.  While the resulting
family of PDE's and numerical schemes are quite general, we give special
attention to their application for regularized inversion problems,
with particular illustrative examples on some popular image
processing applications.
The method is a generalization
of momentum, or accelerated, gradient descent to the PDE setting.
For elliptic problems, the descent equations are a nonlinear damped
wave equation, instead of a diffusion equation, and the acceleration
is realized as an improvement in the CFL condition from
$\Delta t\sim \Delta x^{2}$ (for diffusion) to
$\Delta t\sim \Delta x$ (for wave equations). We work out several
explicit as well as a semi-implicit numerical schemes, together with their
necessary stability constraints, and include recursive update formulations
which allow minimal-effort adaptation of existing gradient descent PDE codes
into the accelerated PDE framework. We explore these schemes more carefully for
a broad class of regularized inversion applications, with special attention to
quadratic, Beltrami, and Total Variation regularization, where the accelerated
PDE takes the form of a nonlinear wave equation. Experimental examples
demonstrate the application of these schemes for image denoising, deblurring,
and inpainting, including comparisons against Primal Dual, Split Bregman, and
ADMM algorithms.
\end{abstract}

\section{Introduction}

Variational problems have found great success, and are widely used,
in image processing for problems such as noisy or blurry image restoration,
image inpainting, image decomposition, and many other problems \cite{aubert2006mathematical}.
Many image processing problems have the form
\begin{equation}
\min_{u}\int_{\Omega}L(x,u,\nabla u)\,dx\label{eq:cvp}
\end{equation}
where $L$ is convex in $\nabla u$\footnote{Nonconvex problems are also widely used, see e.g., \cite{perona1990scale}.}
and the corresponding gradient descent equation
\[
u_{t}+L_{z}(x,u,\nabla u)-\div(\nabla_{p}L(x,u,\nabla u))=0
\]
is a nonlinear diffusion equation, where $L=L(x,z,p)$. Solving \eqref{eq:cvp}
via gradient descent is inefficient, due in large part to the stiff
stability (CFL) condition $\Delta t\leq C\Delta x^{2}$ for diffusion equations.
This has led to the development of more efficient optimization algorithms,
such as primal dual methods \cite{chambolle2011first} and the split
Bregman approach \cite{goldstein2009split} that avoid this numerical
stiffness.\footnote{Primal dual and split Bregman also avoid the non-smoothness of the
$L^{1}$ norm, which is an issue in descent based approaches, which
require some regularization.}

Optimization is also widely used in machine learning, though the types
of optimization problems are (usually) structurally different than
in image processing. As in image processing, for modern large scale
problems in machine learning first order methods based on computing
only the gradient are preferable, since computing and storing the
Hessian is intractable \cite{bottou2010large}. The discrete version
of gradient descent is
\begin{equation}
x_{k+1}=x_{k}-\alpha\nabla f(x_{k}),\label{eq:graddesc}
\end{equation}
where in machine learning the time step $\alpha$ is called the \emph{learning
rate}. While gradient descent is provably convergent for convex problems
\cite{boyd2004convex}, the method can be very slow to converge in
practice.

To address this issue, many versions of \emph{accelerated} gradient
descent have been proposed in the literature, and are widely used
in machine learning \cite{wibisono2016variational}. At some heuristic
level, gradient descent is often slow to converge because the local
descent direction is not reliable on a larger scale, leading to large
steps in poor directions and large corrections in the opposite direction.
Accelerated descent methods incorporate some type of averaging of
past descent directions, which provides a superior descent direction
compared to the local gradient. One of the oldest accelerated methods
is Polyak's heavy ball method \cite{polyak1964some}
\begin{equation}
x_{k+1}=x_{k}-\alpha\nabla f(x_{k})+\beta(x_{k}-x_{k-1}).\label{eq:polyak}
\end{equation}
The term $\beta(x_{k}-x_{k-1})$ acts to average the local descent
direction with the previous direction, and is referred to as \emph{momentum}.
Polyak's heavy ball method was studied in the continuum by Attouch,
Goudou, and Redont~\cite{attouch2000heavy}, and also by Goudou and
Munier \cite{goudou2009gradient}, who call it the \emph{heavy ball
with friction}. In the continuum, Polyak's heavy ball method corresponds
to the equations of motion for a body in a potential field, which
is the second order ODE
\begin{equation}
\ddot{x}+a\dot{x}=-\nabla f(x).\label{eq:polyakode}
\end{equation}
A more recent example of momentum descent is the famous Nesterov accelerated
gradient descent~\cite{nesterov1983method}
\begin{equation}
x_{k+1}=y_{k}-\alpha\nabla f(y_{k}),\ \ y_{k+1}=x_{k+1}+\frac{k-1}{k+2}(x_{k+1}-x_{k}).\label{eq:nes}
\end{equation}
In \cite{nesterov1983method} Nesterov proved a convergence rate of
$O(1/t^{2})$ after $t$ steps for strongly convex problems. This
is provably optimal for first order methods.

The seminal works of Polyak and Nesterov have spawned a whole field
of momentum-based descent methods, and variants of these methods are
widely used in machine learning, such as the training of neural networks
in deep learning \cite{sutskever2013importance,wibisono2016variational}.
The methods are popular for both their superior convergence rates
for convex problems, but also their ability to avoid local minima
in nonconvex problems, which is not fully understood in a rigorous
sense. There has been significant interest recently in understanding
the Nesterov accelerated descent methods. In particular, Su, Boyd
and Candes \cite{su2014differential} recently showed that Nesterov
acceleration is simply a discretization of the second order ODE
\begin{equation}
\ddot{x}+\frac{3}{t}\dot{x}=-\nabla f(x).\label{eq:contNest}
\end{equation}
Other works have since termed this ODE as \emph{continuous time Nesterov}
\cite{ward2017toolkit}. We note the friction coefficient $3/t$ vanishes
as $t\to\infty$, which explains why many implementations of Nesterov
acceleration involve \emph{restarting}, or resetting the time to $t=0$
when the system is underdamped \cite{ward2017toolkit}.

However, it
is the work of Wibisono, Wilson, and Jordan \cite{wibisono2016variational}
that gives the clearest picture of Nesterov acceleration. They show
that virtually all Nesterov accelerated gradient descent methods are
simply discretizations of the ODE equations of motion for a particular
Lagrangian action functional. This endows Nesterov acceleration with
a variational framework, which aids in our understanding, and more
importantly can be easily adapted to other settings.

This was
extended to the partial differential equation (PDE) setting by
Sundaramoorthi and Yezzi in two initial works
\cite{yezzi2017accelerated,sundaramoorthi2018accelerated} (see also
\cite{sundaramoorthi2018accelerated-nips})
where the first set of \emph{Accelerated PDE's} were formulated
both for geometric flows of contours and surfaces (active contours)
as well as for diffeomorphic mappings between images (optical flow).
There are also some acceleration-type
methods that have appeared recently in image processing
\cite{grewenig2016fsi,hafner2016fsi,weickert2016cyclic,bahr2017fast},
however these methods are not derived from a variational framework,
and so they lack energy monotonicity and convergence guarantees.

\subsection{Contributions and Outline of Paper}

In this paper, we extend the class of {\it Accelerated PDE's} formulated in
\cite{yezzi2017accelerated,sundaramoorthi2018accelerated} to
the setting of generic functions over $\R^{n}$, building
on the variational insights pioneered by \cite{wibisono2016variational}.
The method applies to solving general problems
in the calculus of variations. In similar spirit to \eqref{eq:polyakode}
and \eqref{eq:contNest}, the descent equations in PDE acceleration
correspond to a continuous second order flow in time which, for a
broad class of regularized inversion problems to be addressed in Section
\ref{sec:inversion}, take on the specific form of damped nonlinear \emph{wave}
equations rather than the reaction-diffusion equations that arise as their
traditional gradient descent counterparts. Accelerated PDE's can be solved
numerically with simple explicit Euler or semi-implicit Euler schemes
which we develop in Section \ref{sec:numerical}.
Here, acceleration is realized in part through an improvement
in the CFL condition from $\Delta t\sim \Delta x^{2}$ for diffusion equations
(or standard gradient descent), to $\Delta t\sim \Delta x$ for wave equations. In fact, we
will show early on in Section \ref{sec:numerical} that the improvement in
the CFL condition for explicit numerical accelerated PDE schemes
(compared with their gradient descent counterparts)
is a completely general property of accelerated PDE's which applies even
when the wave equation structure explored with more detail in
Section \ref{sec:inversion} does not arise.

In Section \ref{sec:experiments}
we apply the method to quadratic, Beltrami, and Total Variation regularized
problems in image processing including denoising, deblurring, and inpainting,
obtaining results that are comparable
to state of the art methods, such as the split Bregman approach, and ADMM,
and superior to primal dual methods. In a companion paper \cite{calder2018anaccelerated}
we study the PDE acceleration method rigorously and prove a convergence
rate, perform a complexity analysis, and show how to optimally select
the parameters, including the damping coefficient (these results are
summarized in Section \ref{sec:pdeacc}).\\

\subsection{Acknowledgments}

Jeff Calder was supported by NSF-DMS grant 1713691. \\
Anthony Yezzi was supported by NSF-CCF grant 1526848
and ARO W911NF-18-1-0281.

\section{PDE acceleration}

\label{sec:pdeacc}

We now present our PDE acceleration framework, which is based on the
seminal work of \cite{wibisono2016variational,yezzi2017accelerated,sundaramoorthi2018accelerated}
with suitable modifications to image processing problems. We consider
the calculus of variations problem
\[
\min_{u}E[u]:=\int_{\Omega}\Phi(x,\nabla u)+\Psi(x,u)\,dx.
\]
The Euler-Lagrange equation satisfied by minimizers is
\begin{equation}
\nabla E[u]:=\Psi_{z}(x,u)-\div(\nabla\Psi(x,\nabla u))=0,\label{eq:ELEq}
\end{equation}
where $\Phi=\Phi(x,p)$, $\nabla\Phi=\nabla_{p}\Phi$ and $\Psi=\Psi(x,z)$.
We note that the gradient $\nabla E[u]$ satisfies
\begin{equation}
\frac{d}{d\eps}\Big\vert_{\eps=0}E[u+\eps v]=\int_{\Omega}\nabla E[u]\,v\,dx\label{eq:var}
\end{equation}
for all $v$ smooth with compact support, and is often called the
$L^{2}$-gradient due to the presence of the $L^{2}$ inner product
on the right hand side.

We define the action integral
\begin{equation}
J[u]=\int_{t_{0}}^{t_{1}}k(t)\left(\frac{1}{2}\int_{\Omega}\rho u_{t}^{2}\,dx-b(t)E[u]\right)\,dt,\label{eq:action}
\end{equation}
where $k(t)$ and $b(t)$ are time dependent weights, $\rho=\rho(x)$
represents a mass density, and $u=u(x,t)$. Notice the action integral
is the weighted difference between kinetic energy $\tfrac{1}{2}\int\rho u_{t}^{2}\,dx$
and potential energy $E[u]$. The PDE accelerated descent equations
are defined to be the equations of motion in the Lagrangian sense
corresponding to the action $J$. To compute the equations of motion,
we take a variation on $J$ to obtain
\[
0=\frac{d}{d\eps}\Big\vert_{\eps=0}J[u+\eps v]=\int_{t_{0}}^{t_{1}}\int_{\Omega}k(t)\rho u_{t}v_{t}-k(t)b(t)\nabla E[u]\,v\,dx
\]
for smooth $v$ with compact support in $\Omega\times(t_{0},t_{1})$.
Integrating by parts in $t$ we have
\[
0=\int_{t_{0}}^{t_{1}}\int_{\Omega}\left[-\frac{\partial}{\partial t}(k(t)\rho u_{t})-k(t)b(t)\nabla E[u]\right]v\,dx.
\]
Therefore, the PDE accelerated descent equations are
\[
\frac{\partial}{\partial t}(k(t)\rho u_{t})=-k(t)b(t)\nabla E[u].
\]
It is more convenient to define $a(t)=k'(t)/k(t)$ are rewrite the
descent equations as
\begin{equation}
\boxed{u_{tt}+a(t)u_{t}=-b(t)\rho(x)^{-1}\nabla E[u].}\label{eq:PDEacc}
\end{equation}
For image processing problems, there is typically no Dirichlet boundary
condition, so the natural variational boundary condition $\nabla_{p}\Phi(x,\nabla u)\cdot\vb{n}=0$
is imposed on the boundary $\partial\Omega$, where $\vb{n}$ is the
outward normal. Often this reduces to the Neumann condition $\frac{\partial u}{\partial\vb{n}}=0$.

In a companion paper \cite{calder2018anaccelerated} we study the
PDE acceleration descent equation \eqref{eq:PDEacc} rigorously. In
particular, we prove energy monotonicity, and a linear convergence
rate. We summarize the results in Lemma \ref{lem:monotone} and Theorem
\ref{thm:rate}.

\begin{lemma}[Energy monotonicity~\cite{calder2018anaccelerated}]\label{lem:monotone}
Assume $a(t),b(t)\geq0$ and let $u$ satisfy \eqref{eq:PDEacc}.
Suppose either $u(x,t)=g(x)$ or $\nabla\Phi(x,\nabla u)\cdot\vb{n}=0$
on $\partial\Omega$. Then
\begin{equation}
\frac{d}{dt}\left(K[u]+b(t)E[u]\right)=-2a(t)K[u]+b'(t)E[u],\label{eq:monotone}
\end{equation}
where $K[u]=\frac{1}{2}\int_{\Omega}\rho u_{t}^{2}\,dx$. In particular,
total energy is always decreasing provided $b'(t)\leq0$ and $E[u]\geq0$.
\end{lemma}
\begin{theorem}[Convergence rate~\cite{calder2018anaccelerated}]\label{thm:rate}
Let $u$ satisfy \eqref{eq:PDEacc} and let $u^{*}$ be a solution
of $\nabla E[u^{*}]=0$ in $\Omega$. Assume $\Phi$ is uniformly
convex in $\nabla u$, $\Psi$ is convex, and $\Psi_{zz}$ is bounded
above, $u=u^{*}$ on $\partial\Omega$, $a(t)=a>0$ is constant and
$b(t)\equiv1$ and $\rho\equiv1$. Then there exists $C,\beta>0$
such that
\begin{equation}
\|u-u^{*}\|_{H^{1}(\Omega)}^{2}\leq C\exp\left(-\beta t\right).\label{eq:rate}
\end{equation}
\end{theorem}

We mention that the same convergence rate \eqref{eq:rate} holds for
gradient descent
\[
u_{t}=-\nabla E[u]
\]
under the same conditions on $E$. The difference is that gradient
descent is a diffusion equation, which requires a times step of $\Delta t\sim \Delta x^{2}$
for stability, while PDE acceleration \eqref{eq:PDEacc} is a wave
equation which allows much larger time steps $\Delta t\sim \Delta x$. Thus, the
acceleration is realized as a relaxation in the CFL condition.

While Theorem \ref{thm:rate} provides a convergence rate, it does
not give advice on how to select the damping coefficient $a>0$. It
was shown in \cite{calder2018anaccelerated} how to optimally select
the damping coefficient in the linear setting, and we find this choice
is useful for nonlinear problems as well. For convenience, we recall
the results from \cite{calder2018anaccelerated}, which apply to the
linear PDE acceleration equation
\begin{equation}
u_{tt}+au_{t}+Lu+\lambda u=f\ \ \text{ in }\Omega\times(0,\infty),\label{eq:Lpmotion}
\end{equation}
where $L$ is a linear second order elliptic operator. A Fourier analysis
\cite{calder2018anaccelerated} leads to the optimal choice
\begin{equation}
a=2\sqrt{\lambda_{1}+\lambda},\label{eq:aopt}
\end{equation}
where $\lambda_{1}$ is the first Dirichlet eigenvalue of $L$ (or
for the Neumann problem, the first eigenvalue corresponding to a nontrivial
eigenfunction), and the optimal convergence rate
\begin{equation}
|u(x,t)-u^{*}(x)|\leq C\exp\left(-at\right).\label{eq:error}
\end{equation}
Notice that if $L$ is degenerate elliptic, so $\lambda_{1}=0$, which
roughly corresponds to a non-strongly convex optimization problem,
the method still converges when $\lambda>0$. That is, the presence
of a fidelity term in the image processing problem enables, and accelerates,
convergence. This suggests why the algorithm is successful even for
TV restoration, which is not strongly convex but has a fidelity.

\section{Numerical Schemes for Accelerated PDE's \label{sec:numerical}}

We now describe various time discretization strategies for the generic accelerated
PDE 
\begin{align}
u_{tt}+au_{t} & =-\nabla E[u,u_{x},u_{xx},\cdots]\quad\mbox{accelerated PDE}\label{eq:PDE-accelerated}\\
u_{t} & =-\nabla E[u,u_{x},u_{xx},\cdots]\quad\mbox{gradient descent PDE}\label{eq:PDE-gradient}
\end{align}
(alongside related discretizations of the generic gradient descent PDE for comparison)
where \eqref{eq:PDE-accelerated} represents the unit density ($\rho=1$)
and unit energy scaled ($b=1$) case of \eqref{eq:PDEacc}. A key advantage
of accelerated PDE schemes for regularized inversion problems, which
we explore subsequently in Section \ref{sec:inversion}, is that in
typical cases where the gradient descent PDE \eqref{eq:PDE-gradient}
takes the form of a linear or nonlinear reaction-diffusion equation,
the matching accelerated PDE \eqref{eq:PDE-accelerated} takes the
form of a linear or nonlinear wave equation, whose explicit time discretization
permits a much larger stable time step than the explicit discretization
of \eqref{eq:PDE-gradient}. Therefore, due to their simplicity of implementation,
as well as their immediately parallelizable structure, we will restrict
our discussion to explicit update schemes and to the semi-implicit
Euler scheme whose two-part update consists of partial updates which
are both explicit in nature.

\subsection{Explicit Forward-Euler for Gradient Descent PDE's}

We start by considering the explicit forward Euler discretization
of the continuous gradient descent PDE \eqref{eq:PDE-gradient}. Using
a forward difference in time to approximate the time derivative on
the left hand side we obtain
\[
\frac{u(x,t+\Delta t)-u(x,t)}{\Delta t}=-\nabla E
\]
This leads to the following simple discrete iteration
\begin{eqnarray}
\Delta u^{n}(x) & = & -\Delta t\,\nabla E^{n}\label{eq:FDE-grad}\\
u^{n+1}(x) & = & u^{n}(x)+\Delta u^{n}(x)\nonumber
\end{eqnarray}
where $u^{n}(x)\doteq u(x,n\Delta t)$ denotes the current iterate,
$\Delta u^{n}\doteq u(x,n\Delta t+\Delta t)-u(x,n\Delta t)$
the increment to be applied, $u^{n+1}(x)\doteq u\left(x,(n+1)\Delta t\right)$
the new iterate, and $\nabla E^{n}(x)\doteq\nabla E(x,n\Delta t)$
the discrete approximation of the gradient computed at step
$n$.

In most cases, stability considerations require an upper bound
on the time step $\Delta t$ (the CFL condition)
dependent upon the discretization of $\nabla E^{n}$.
Often this upper bound for stable time steps is computed using Von
Neumann analysis by linearizing $\nabla E^{n}$ in \eqref{eq:FDE-grad}
and taking a Discrete Fourier Transform (DFT) on both sides of the
homogeneous part to obtain.
\[
U^{n+1}(\omega)-U^{n}(\omega)=-\Delta t\left[z(\omega)U^{n}(\omega)\right]
\]
Such a structure often arises when $\nabla E^{n}$ is computed explicitly
using only the values of $u^{n}$. In such cases, its linearization
will consist of a combination of $u^{n}$ values whose DFT can be
written in the form $z(\omega)\,U^{n}(\omega)$ where $U^{n}(\omega)$
denotes the DFT of $u^{n}$. We will refer to $z(\omega)$ as the\emph{
gradient amplifier}.\footnote{
  A discrete version of what is often called the \emph{symbol} of the
  underlying linear differential operator that is being approximated.
}

This leads to the following update
\begin{eqnarray}
U^{n+1}(\omega) & = & \underbrace{\left(1-\Delta t\,z(\omega)\right)}_{\xi(\omega)}U^{n}(\omega)\nonumber \\
\mbox{gradient amplifier }z(\omega) & \doteq & \frac{\mbox{DFT}\left(\mbox{linearized\_homogeneous\_part\_of}\left(\nabla E^{n}\right)\right)}{\mbox{DFT}\left(u^{n}\right)}\label{eq:grad_amplifier}
\end{eqnarray}
which will be stable as long as the overall update amplification factor
$\xi(\omega)$ does not have complex amplitude exceeding unity for
any frequency $\omega$. This condition can be expressed as
\begin{align*}
\xi(\omega)\xi^{*}(\omega) & =\left(1-\Delta t\,z(\omega)\right)\left(1-\Delta t\,z^{*}(\omega)\right)\le1\\
\Delta t & \le\frac{z(\omega)+z^{*}(\omega)}{z(\omega)z^{*}(\omega)}=\frac{1}{z(\omega)}+\frac{1}{z^{*}(\omega)}=2\Re\left(\frac{1}{z(\omega)}\right)
\end{align*}
For elliptic operators, which are common in regularized optimization in image
processing, the gradient amplifier is real and non-negative:
$z(\omega)\ge 0$. In such cases the stability
constraint takes the form of the following CFL condition
\begin{equation}
\Delta t\le\frac{2}{z_{\max}}\label{eq:CFL-grad}
\end{equation}
where $z_{\max}\doteq\max_{\omega}z(\omega)$.

\subsection{Fully Explicit Schemes for Accelerated PDE's}

We now turn our attention to the explicit discretizations of the accelerated
PDE \eqref{eq:PDE-accelerated}. We will consider both first and second
order approximations of the time derivatives and will exploit the
following lemma in the Von Neumann stability analysis for each of
these choices.

\paragraph{Root Amplitude Lemma}

\emph{Given a quadratic equation $A\xi^{2}+B\xi+C=0$ with real coefficients
($A\ne0)$, its roots will satisfy $|\xi|\le1$ if and only if $\frac{\left|B\right|}{|A|}-1\le\frac{C}{A}\;\le1$
(or equivalently $A>C$ and $A+C>|B|$ for positive $A$). }
\begin{quote}
\emph{\small{}Proof.}{\small{} We first prove the result in the special
case that $A=1$ and $B\ge0$, in which case the roots are $\xi=-\frac{B}{2}\pm\frac{1}{2}\sqrt{B^{2}-4C}$
and claim that $|\xi|\le1$ if and only if
\[
B-1\le C\le1
\]
If the roots are imaginary, then both have complex amplitude
$|\xi|^{2}=C>0$ which makes the right hypothesis both necessary and
sufficient. The left hypothesis automatically follows since $C>\frac{B^{2}}{4}\ge B-1$
(the first part for the roots to be imaginary and the second part
equivalent to $(B-2)^{2}\ge0$). In the case of real roots, we want
the larger magnitude root to satisfy $|\xi|=\frac{B}{2}+\frac{1}{2}\sqrt{B^{2}-4C}\le1$,
which can be expressed as $\sqrt{B^{2}-4C}\le2-B$. This immediately
yields $B\le2$ as a necessary condition to keep the right side positive.
Under this condition, we can square both sides and simplify to obtain
the left hypothesis as necessary and sufficient. The right hypothesis
automatically follows since $C<\frac{B^{2}}{4}<1$ (the first part
for the roots to be real and the second part based on our condition).
Combining the hypotheses yields $B-1\le1$ which satisfies the necessary
condition, thus completing the special case proof. The general case
follows since the roots of $A\xi^{2}+B\xi+C$ have the same magnitude
as the roots of $\xi^{2}+\frac{|B|}{|A|}\xi+\frac{C}{A}$.}{\small\par}
\end{quote}

\subsubsection{Second order in time scheme}

Using central difference approximations for both time derivatives
gives a second order discretization in time
\begin{align*}
\frac{u(x,t+\Delta t)-2u(x,t)+u(x,t-\Delta t)}{\Delta t^{2}}+a\frac{u(x,t+\Delta t)-u(x,t-\Delta t)}{2\Delta t} & =-\nabla E(x,t)
\end{align*}
which leads to the following update.
\begin{equation}
u^{n+1}(x)=\frac{2u^{n}(x)-\left(1-\frac{a\Delta t}{2}\right)u^{n-1}(x)-\Delta t^{2}\nabla E^{n}(x)}{1+\frac{a\Delta t}{2}}\label{eq:FDE-central}
\end{equation}
Applying the DFT to the linearized homogeneous part of the update
scheme \eqref{eq:FDE-central} yields
\begin{align*}
U^{n+1}(\omega) & =\frac{\left(2-\Delta t^{2}\,z(\omega)\right)U^{n}(\omega)-\left(1-\frac{a\Delta t}{2}\right)U^{n-1}(\omega)}{1+\frac{a\Delta t}{2}}
\end{align*}
where $z(\omega)$ denotes the gradient amplifier \eqref{eq:grad_amplifier}.
If we substitute $U^{n\pm m}=\xi^{\pm m}U^{n}$, where $\xi(\omega)$
denotes the overall update amplification factor, then we obtain the
quadratic equation
\begin{align*}
\underbrace{\left(1+\frac{a\Delta t}{2}\right)}_{A}\xi^{2} & +\underbrace{\left(\Delta t^{2}\,z(\omega)-2\right)}_{B}\xi+\underbrace{\left(1-\frac{a\Delta t}{2}\right)}_{C}=0
\end{align*}

In the case of real $z(\omega)$, we may exploit the Root Amplitude
Lemma to check the stability condition $|\xi(\omega)|\le1$. The first
condition $A\ge C$ of the lemma (for positive $A$) is satisfied
since $1+\frac{a\Delta t}{2}\ge1-\frac{a\Delta t}{2}$ for all positive
$a$ and $\Delta t$, and so we use the second condition $A+C\ge|B|$
to obtain the stability condition $2\ge\left|2-\Delta t^{2}z(\omega)\right|$,
which may be rewritten as $0\le\Delta t^{2}z(\omega)\le4$. In the
case where $z(\omega)\ge0$ we automatically satisfy the left hand
inequality for all $\omega$, leaving us with
\begin{equation}
\Delta t\le\frac{2}{\sqrt{z_{\max}}}\label{eq:CFL-central}
\end{equation}

\subsubsection{First order in time schemes }

Continuing to use a central difference for the second derivative but
only a one sided difference (forward or backward) for the first derivative
in time, yields two alternative first order time schemes.

\paragraph{Forward difference}

\begin{align*}
\frac{u(x,t+\Delta t)-2u(x,t)+u(x,t-\Delta t)}{\Delta t^{2}}+a\frac{u(x,t+\Delta t)-u(x,t)}{\Delta t} & =-\nabla E(x,t)
\end{align*}
The forward difference discretization yields the update.
\begin{equation}
u^{n+1}(x)=\frac{\left(2+a\Delta t\right)u^{n}(x)-u^{n-1}(x)-\Delta t^{2}\nabla E^{n}(x)}{1+a\Delta t}\label{eq:FDE-forward}
\end{equation}
Von Neumann analysis applied to the linearized homogeneous part of
\eqref{eq:FDE-central} yields the following quadratic equation for
the update amplification factor $\xi(\omega)$.
\begin{align*}
\underbrace{\left(1+a\Delta t\right)}_{A}\xi^{2}(\omega) & +\underbrace{\left(\Delta t^{2}\,z(\omega)-2-a\Delta t\right)}_{B}\xi(\omega)+\underbrace{1}_{C}=0
\end{align*}
Since $1+a\Delta t>1$ for all positive $a$ and $\Delta t$, the
first condition $A\ge C$ of the root amplitude lemma (for positive
A) is always satisfied. We may therefore restrict out attention to
the second condition $A+C\ge|B|$, assuming real $z(\omega)$, to
determine whether $|\xi(\omega)|\le1$. This gives the condition $(1+a\Delta t)+1\ge\left|\Delta t^{2}z(\omega)-(2+a\Delta t)\right|$
which is equivalent to $0\le\Delta t^{2}z(\omega)\le2(2+a\Delta t)$.
In the case where $z(\omega)\ge0$ we automatically satisfy the left
hand inequality for all $\omega$, which leaves us with $z(\omega)\,\Delta t^{2}-2a\Delta t-4\le0$.
Plugging in the extreme case $z_{\max}$ and restricting $\Delta t$
to lie below the positive root in order to keep the quadratic expression
on the left negative yields

\begin{equation}
\Delta t\le\sqrt{\frac{4}{z_{\max}}+\left(\frac{a}{z_{\max}}\right)^{2}}+\frac{a}{z_{\max}}\label{eq:CFL-forward}
\end{equation}
Notice that the CFL condition \eqref{eq:CFL-central} for the second
order (central difference) scheme is sufficient but not necessary.
If, however, we wish to obtain a condition independent of the damping
$a$, then minimizing the upper bound with respect to $a$ (by plugging
in $a=0$) recovers this prior second order CFL condition.

\paragraph{Backward difference}

\begin{align*}
\frac{u(x,t+\Delta t)-2u(x,t)+u(x,t-\Delta t)}{\Delta t^{2}}+a\frac{u(x,t)-u(x,t-\Delta t)}{\Delta t} & =-\nabla E(x,t)
\end{align*}
The backward difference discretization yields the update.
\begin{equation}
u^{n+1}(x)=\left(2-a\Delta t\right)u^{n}(x)-\left(1-a\Delta t\right)u^{n-1}(x)-\Delta t^{2}\nabla E^{n}(x)\label{eq:FDE-backward}
\end{equation}
Similar analysis yields the following quadratic equation for the amplification
factor $\xi(\omega)$.
\begin{align*}
\xi^{2}(\omega)-\left(2-a\Delta t-\Delta t^{2}z(\omega)\right)\xi(\omega)+(1-a\Delta t) & =0
\end{align*}
The first condition $A\ge C$ of the lemma (for positive $A$) is
always satisfied ($1\ge1-a\Delta t)$ for all positive values of $a$
and $\Delta t$. The second condition $A+C\ge|B|$, assuming real
$z(\omega)\ge0$, of the lemma, can be expressed as $\Delta t^{2}z(\omega)+2a\Delta t-4\le0.$
Plugging in the extreme case $z_{\max}$ and restricting $\Delta t$
to lie below the positive root in order to keep the quadratic expression
on the left negative gives the following CFL condition.

\begin{equation}
\Delta t\le\sqrt{\frac{4}{z_{\max}}+\left(\frac{a}{z_{\max}}\right)^{2}}-\frac{a}{z_{\max}}\label{eq:CFL-backward}
\end{equation}
Notice that the CFL condition \eqref{eq:CFL-central} for the central
difference scheme is necessary (easily seen by applying the triangle
inequality) but no longer sufficient as in the forward difference
case. Furthermore, the constraint becomes increasingly restrictive
as the damping coefficient $a$ increases, making it impossible to
formulate a sufficient damping-independent stability constraint. We
will therefore give no further consideration to this scheme.

\subsection{Recursive Increments and Properties of Explicit Schemes}

For greater convenience in implementation, especially when upgrading
existing gradient descent routines structured according to \eqref{eq:FDE-grad}
with one array to store the evolving iterate $u^{n}$ and another
for its increment $\Delta u^{n}$, the explicit accelerated PDE discretizations
can be expressed in terms of recursively defined increments as follows\footnote{For completeness, the first order backward difference scheme can also
be written recursively in the form $\Delta u^{n}=\left(1-a\Delta t\right)\Delta u^{n-1}-\Delta t^{2}\nabla E^{n}$.}
\begin{align}
\mbox{\footnotesize gradient descent: }\Delta u^{n} & =\qquad\qquad\qquad\qquad-\Delta t\,\nabla E^{n},\qquad\Delta t\le\frac{2}{z_{\max}}\label{eq:update-gradient}\\
\mbox{\footnotesize 1-order accelerated: }\Delta u^{n} & =\frac{1}{1+a\Delta t}\Delta u^{n-1}-\frac{\Delta t^{2}}{1+a\Delta t}\nabla E^{n},\;\Delta t\le\sqrt{\frac{4}{z_{\max}}+\left(\frac{a}{z_{\max}}\right)^{2}}+\frac{a}{z_{\max}}\label{eq:update-forward}\\
\mbox{\footnotesize 2-order accelerated: }\Delta u^{n} & =\frac{2-a\Delta t}{2+a\Delta t}\Delta u^{n-1}-\frac{2\Delta t^{2}}{2+a\Delta t}\nabla E^{n},\;\Delta t\le\frac{2}{\sqrt{z_{\max}}}\label{eq:update-central}
\end{align}
where $\Delta u^{n-1}$ denotes the previously increment (kept in
just one more added array). Here we see more directly the traditional
momentum style structure (i.e. heavy-ball) in that the next increment
$\Delta u^{n}$ is expressed as a weighted combination of the gradient
$\nabla E^{n}$ and the previous increment $\Delta u^{n-1}$. Recursion
\eqref{eq:update-forward} is equivalent to the first order in time,
explicit update \eqref{eq:FDE-forward} using forward differences
while the recursion \eqref{eq:update-central} is equivalent to the
second order in time, explicit update \eqref{eq:FDE-central} using
central differences, and as such they must adhere to the same CFL
conditions \eqref{eq:CFL-forward} and \eqref{eq:CFL-central} derived
earlier for these corresponding schemes.

\subsubsection{The first order scheme as a sub-case of the second order scheme}

For any choice of damping $\alpha_{1}$ and time step $\Delta t_{1}$
parameters used in the first order scheme (denoted by subscript 1),
we may obtain equivalent update iterations by substituting the following
change of parameters into the second order scheme (denoted by subscript
2).
\begin{equation}
\Delta t_{2}=\frac{\Delta t_{1}}{\sqrt{1+\frac{a_{1}\Delta t_{1}}{2}}}\quad\mbox{and}\quad a_{2}=\frac{a_{1}}{\sqrt{1+\frac{a_{1}\Delta t_{1}}{2}}}\label{eq:first-to-second}
\end{equation}
This is easily shown by algebraic simplification of the second order
update (and stability condition) after applying the change of parameters.
The simplified result will yield the first order scheme (and stability
condition) in the original damping and time step parameters. In short,
the first order scheme is always equivalent to the second order scheme
with a reduced time step and damping via the contraction factor $\frac{1}{\sqrt{1+\frac{a_{1}\Delta t_{1}}{2}}}<1$.

A particular special case of this equivalency arises in considering
the maximal stable time step for both schemes. For a fixed choice
of damping $a$, the first order scheme appears to allow a more generous
upper bound than the second order scheme. However, there is no effective
difference when substituting \eqref{eq:first-to-second} into the
second order scheme. Although the upper bound on the time step is
smaller, the contracted time step is also smaller, such that the maximum
stable time step in the first order scheme rescales exactly to the
maximum stable time step in the second order scheme. Thus, so long
as the damping is also contracted according to \eqref{eq:first-to-second},
the first order scheme implemented with its maximum stable time step
is equivalent to the second order scheme implemented with its maximum
stable time step.

We may also consider the backwards version of the change of parameters
\eqref{eq:first-to-second} in order to map the second order scheme
into the first order scheme. In this case, using parameters $\alpha_{2}$
and $\Delta t_{2}$ in the second order scheme is equivalent to applying
the following change of parameters to the first order scheme.
\begin{equation}
\Delta t_{1}=\frac{\Delta t_{2}}{\sqrt{1-\frac{a_{2}\Delta t_{2}}{2}}}\quad\mbox{and}\quad a_{1}=\frac{a_{2}}{\sqrt{1-\frac{a_{2}\Delta t_{2}}{2}}}\label{eq:second-to-first}
\end{equation}
However, this backwards mapping only applies when the second order
parameters satisfy $a_{2}\Delta t_{2}<2$. When this condition is
satisfied, we can show by direct substitution and algebraic simplification,
that the second order scheme (and stability condition) is equivalent
to the first order scheme (and stability condition) with an amplified
time step and damping coefficient via the amplification factor $\frac{1}{\sqrt{1-\frac{a_{2}\Delta t_{2}}{2}}}>1$.
Assuming the same condition is satisfied, the second order scheme
implemented with its maximum stable time step is equivalent to the
first order scheme implemented with its maximum stable time step after
boosting the damping parameter according to \eqref{eq:second-to-first}.

\subsubsection{Critical damping in the second order scheme (gradient descent)}

Unlike the forward mapping of the first order into the second order
discrete scheme, which is always possible for any choice of first
order discrete parameters $a_{1}$ and $\Delta t_{1}$, the backwards
mapping is not possible for certain choices of the second order discrete
parameters, namely for $a_{2}\Delta t_{2}\ge2$ where the backwards
amplification factor $\frac{1}{\sqrt{1-\frac{a_{2}\Delta t_{2}}{2}}}$
is undefined. While $\Delta t_{2}$ is upper-bounded by the second
order scheme's stability constraint, there is no such upper-bound
imposed on $\alpha_{2}$ since the stability constraint is independent
of $\alpha_{2}$. As such for any stable, nonzero, second-order discrete
time step $\Delta t_{2}$, we may always choose the second order discrete
damping coefficient $a_{2}$ high enough to enter into this parameter
regime where $a_{2}\Delta t_{2}\ge2$. In this case the second order
scheme will exhibit behavior that is no longer reproducible by the
first order scheme.

It is interesting to consider what happens at the transition point
when $a_{2}\Delta t_{2}=2$. It is immediately seen, by plugging this
into \eqref{eq:update-central}, that the second order scheme becomes
identical to the discrete gradient descent scheme \eqref{eq:update-gradient}
with an effective gradient descent time step of $\Delta t=\frac{1}{2}\Delta t_{2}^{2}$
at this transition point (and if the second order accelerated time
step $\Delta t_{2}$ was chosen to be the maximum stable step size
of $2/\sqrt{z_{\max}}$, the effective gradient descent time step
$\Delta t$ will also be the maximum stable gradient descent step
size of $2/z_{\max}$). If we fix the second order step size $\Delta t_{2}$
and approach the transition point $a_{2}=2/\Delta t_{2}$ from below,
where an equivalent first order damping coefficient $a_{1}$ can be
obtained via \eqref{eq:second-to-first}, then we see that the damping
in the matching first order scheme becomes infinite as the damping
in the second order scheme approaches this critical value. This constitutes
a discrete analog of the continuum property that the continuous gradient
descent PDE \eqref{eq:PDE-gradient} arises as the infinite frictional
limit of the continuous accelerated descent PDE \eqref{eq:PDE-accelerated}.

If we want a damping value $\alpha_{2}$ in the second order scheme
that will always keep us below this transition point for all choices
of stable time step, then we must satisfy the inequality $a_{2}\Delta t_{2}<2$
for the maximum stable step size of $2/\sqrt{z_{\max}}$. This leads
to the following upper bound for the second order damping coefficient.
\[
a_{2}<\sqrt{z_{\max}}
\]
Namely, the square of damping factor should be strictly less than
the gradient amplifier.

\subsubsection{Over-damping in the second order scheme (gradient descent with resistance)}

Noting that gradient descent arises in both schemes (although only
in the limiting sense for the first order scheme) at the transition
point when $a_{2}\Delta t_{2}=2$, and that both schemes offer equivalent
discretizations of accelerated descent according to the rescalings
\eqref{eq:first-to-second} and \eqref{eq:second-to-first} below
this transition point, it is now interesting to consider what happens
above this transition point in the second order scheme. If we choose
$a_{2}>\sqrt{z_{\max}}$ then there will be stable time step choices
for $\Delta t_{2}$ that will bring us beyond this transition.

In the case $a_{2}\Delta t_{2}>2$, the second order update, in its
recursive form \eqref{eq:update-central} becomes a weighted combination
of a step in the negative gradient direction as well as a \emph{backward}
step in the previous update direction. As such, the combined step
can be interpreted as \emph{partially undoing} the previous step,
thereby slowing down the descent process. If we take the limiting
case as the second order damping coefficient $a_{2}$ becomes infinite
(keeping the same fixed time step $\Delta t_{2}$), the stability
of the scheme will not be affected, but the new update will fully
subtract the previous update, thereby returning to the previous state
before applying the new gradient step. Furthermore, after subtracting
the previous update the amount of movement along the new gradient
step will be zero for infinite $a_{2}$. This can be seen by noting
that the weight on the previous update in \eqref{eq:update-central}
approaches -1 from above and that the weight on the gradient approaches
0 from below as $a_{2}\to\infty$. Therefore, in the limit, even if
we initialize the recursion \eqref{eq:update-central} with a non-zero
starting update $\Delta u^{0}$ (the discrete analog of an initial
velocity), the effect will still be to remain motionless at the initial
condition $u^{0}$.

This leads to the interpretation of the over-damped case as a \emph{resisted}
version of gradient descent for any finite $\alpha_{2}>2/\Delta t_{2}$,
since we start with a gradient step in the first update, then partially
undo it before taking a new gradient step in the second update, which
is then partially undone before taking a new gradient step in the
third update, and so on. Since the fraction of each gradient step
which gets subtracted in the subsequent step remains fixed, rather
than accumulating, we do not refer to this as deceleration but rather
as \emph{resistance}, which impedes the normal progress of gradient
descent by a constant factor. As $a_{2}$ increases, resistance increases,
further slowing the progress of gradient descent, while completely
halting it in the limit as $a_{2}\to\infty$.

\subsection{Semi-Implicit Schemes}

We may use semi-implicit Euler style discretizations of \eqref{eq:PDE-accelerated}
to obtain systems which more closely resemble the classic two-part
Nesterov recursion. We may do this with any of the fully explicit
schemes \eqref{eq:FDE-central}, \eqref{eq:FDE-forward}, or \eqref{eq:FDE-backward}
by replacing the explicit discretization $\nabla E^{n}$ of the gradient
with a ``predicted estimate'' $\widehat{\nabla E}^{n+1}$of its
implicit discretization $\nabla E^{n+1}$. This estimate is obtained
by applying the same discretization of $\nabla E$ used in approximating
$\nabla E^{n}\approx\nabla E(u^{n})$ to a partial update $v^{n}$
for the ``look ahead'' approximation $\widehat{\nabla E}^{n+1}\approx\nabla E(v^{n})$.
The partial update $v^{n}$ is obtained beforehand via the fully explicit
update without the gradient term (i.e. by treating $\nabla E^{n}$
as if it were zero). Using this strategy with the second order in
time scheme \eqref{eq:FDE-central} yields the following two-step
update, where the first and second steps, in isolation, both have
a fully explicit structure.

\begin{equation}
v^{n}=u^{n}+\frac{2-a\Delta t}{2+a\Delta t}\Delta u^{n-1}\quad\mbox{then}\quad u^{n+1}=v^{n}-\frac{2\Delta t^{2}}{2+a\Delta t}\,\underbrace{\nabla E\left(v^{n}\right)}_{\approx\widehat{\nabla E}^{n+1}}\label{eq:FDE-semi}
\end{equation}

Von Neumann analysis can be employed to analyze the stability of this
scheme according the following update relationships between the DFT
sequences $U^{n}$, $V^{n}$, and $U^{n+1}$ (transforms of $u^{n}$,
$v^{n}$, and $u^{n+1}$ respectively) where $z(w)$ represents the
gradient amplifier \eqref{eq:grad_amplifier} associated with the
linearization of $\nabla E^{n}$ (and therefore also with the linearization
of $\widehat{\nabla E}^{n+1}$).
\begin{align*}
V^{n} & =U^{n}+\frac{2-a\Delta t}{2+a\Delta t}\left(U^{n}-U^{n-1}\right)=\frac{4}{2+a\Delta t}U^{n}-\frac{2-a\Delta t}{2+a\Delta t}U^{n-1}\\
U^{n+1} & =V^{n}-\frac{2\Delta t^{2}}{2+a\Delta t}z(\omega)V^{n}=\left(1-\frac{2\Delta t^{2}}{2+a\Delta t}z(\omega)\right)V^{n}
\end{align*}
If we substitute the first expression into the second, followed by
substitutions $U^{n\pm m}=\xi^{\pm m}U^{n}$, then we obtain the quadratic
equation
\[
\underbrace{\left(2+a\Delta t\right)^{2}}_{A}\xi^{2}-\underbrace{4\left(2+a\Delta t-2\Delta t^{2}z(\omega)\right)}_{-B}\xi+\underbrace{\left(2-a\Delta t\right)\left(2+a\Delta t-2\Delta t^{2}z(\omega)\right)}_{C}=0
\]
for the overall combined update amplification factor $\xi(\omega)$.
We may use the Root Amplitude Lemma to check the stability criterion
$|\xi(\omega)|\le1$.

\paragraph{First stability condition: $A\ge C$}

The first condition from the lemma (for positive $A$) can be expressed
in quadratic form as $az\Delta t^{2}-\left(a^{2}+2z\right)\Delta t-2a\le0$
which will be satisfied between its positive and negative roots. Restricting
our interest to only positive values of $\Delta t$ therefore yields
the constraint
\[
\Delta t\le\underbrace{\frac{a^{2}+2z+\sqrt{\left(a^{2}+2z\right)^{2}+8a^{2}z}}{2az}}_{g(a,z)>0}\quad\text{where }\frac{\partial g}{\partial a}=\underbrace{\left(a^{2}-2z\right)}_{\mbox{toggles}}\left(\underbrace{\frac{a^{2}+2z+\sqrt{\left(a^{2}+2z\right)^{2}+8a^{2}z}}{2a^{2}z\sqrt{\left(a^{2}+2z\right)^{2}+8a^{2}z}}}_{\mbox{always positive}}\right)
\]
To satisfy this independently of $a$, we examine the partial derivative
of the upper bound $g(a,z)$ with respect to $a$ see that it starts
out negative for $a^{2}<2z$ then turns positive for $a^{2}>2z$.
The minimum upper bound is therefore attained when $z(\omega)=z_{\max}$
and $a^{2}=2z_{\max}$ yielding
\[
\Delta t\le\frac{2+\sqrt{2}}{\sqrt{z_{\max}}}
\]
While this upper bound is more generous than \eqref{eq:CFL-central}
for the fully explicit scheme, it only satisfies the first of the
two stability conditions in the Bounded Root Lemma. We now proceed
to the second condition which will be more restrictive.

\paragraph{Second stability condition: $A+C\ge|B|$}

The second condition from the lemma (for positive $A$) can be expressed
as

\[
2+a\Delta t-\Delta t^{2}z(\omega)+\frac{1}{2}a\Delta t^{3}z(\omega)\ge\left|2+a\Delta t-2\Delta t^{2}z(\omega)\right|
\]
For small enough time steps, $2+a\Delta t-2\Delta t^{2}z$ is positive,
the absolute value signs can be removed, and the inequality holds.
For larger time steps $2+a\Delta t-2\Delta t^{2}z$ becomes and the
inequality can be rearranged into the following cubic form
\[
az(\omega)\Delta t^{3}-6z(\omega)\Delta t^{2}+4a\Delta t+8\ge0
\]
Minimizing on the left with the case $a=0$ and $z(\omega)=z_{\max}$
yields a stricter, and therefore sufficient, stand-alone stability
condition
\begin{equation}
\Delta t\le\frac{2}{\sqrt{3z_{\max}}}\label{eq:CFL-semi}
\end{equation}
Note that this upper bound is smaller, by a factor of $\sqrt{3}$,
than the maximum stable time step \eqref{eq:CFL-central} for the
corresponding fully explicit scheme \eqref{eq:FDE-central} or for
its recursive equivalent \eqref{eq:update-central}.

\section{Regularized Inversion via Accelerated PDEs\label{sec:inversion}}

Here we consider a very general class of variational regularized inversion
problems in the accelerated PDE framework. In particular, we assume
energy functions with the form
\[
E(u)=\int_{\Omega}\underbrace{f\left(|{\cal K}u-g|\right)}_{\mbox{fidelity}}+\underbrace{r(\|\nabla u\|)}_{\mbox{regularity}}\,dx,\qquad\mbox{with }\dot{f},\dot{r},\ddot{r}>0
\]
where $f$ is a monotonically increasing penalty on the residual error
between data measurements $g$ and a forward in the form of linear
operator ${\cal K}$ applied to the reconstructed signal $u$, while
$r$ is a monotonically increasing penalty on the gradient of the
reconstruction.

\subsection{General Case (nonlinear wave equation)\label{sec:general}}

The continuum gradient of $E$ has the form
\begin{align*}
\nabla E(u) & =\frac{\dot{f}\left(|{\cal K}u-g|\right)}{|{\cal K}u-g|}{\cal K}^{*}\left({\cal K}u-g\right)-\dot{r}(\|\nabla u\|)\,\nabla\cdot\left(\frac{\nabla u}{\|\nabla u\|}\right)-\,\ddot{r}(\|\nabla u\|)\,\frac{\nabla u^{T}\nabla^{2}u\,\nabla u}{\|\nabla u\|^{2}}\\
 & =\underbrace{\frac{\dot{f}\left(|{\cal K}u-g|\right)}{|{\cal K}u-g|}}_{\lambda(u,x)>0}{\cal K}^{*}\left({\cal K}u-g\right)-\underbrace{\frac{\dot{r}(\|\nabla u\|)}{\|\nabla u\|}}_{c(\nabla u)>0}\left(\nabla\cdot\nabla u-u_{\eta\eta}\right)-\underbrace{\ddot{r}(\|\nabla u\|)}_{d(\nabla u)>0}u_{\eta\eta}
\end{align*}
where ${\cal K}^{*}$ denotes the adjoint of the forward operator
${\cal K}$, and where $\eta\doteq\frac{\nabla u}{\|\nabla u\|}$
denotes the unit vector along the gradient direction of $u$. This
gives rise to the following class of accelerated flows which take
the form of a nonlinear wave equation.
\begin{equation}
u_{tt}-\,c(\nabla u)\,\left(\nabla\cdot\nabla u-u_{\eta\eta}\right)\,-\,d(\nabla u)\;u_{\eta\eta}+au_{t}=\lambda(u,x)\;{\cal K}^{*}\left(g-{\cal K}u\right)\,
\end{equation}

If, purely for the sake of understanding stability, we model the short
time behavior of any of the presented discrete update schemes in the
neighborhood of a particular spatial point $x$, by treating $\lambda$,
$c$, and $d$ as locally constant, and by representing the forward
model linear operator ${\cal K}$ as a real convolution kernel $K$
with adjoint $K^{T}$, then $\nabla E$ can be approximated near
$x$ by the following linear expression
\begin{equation}
\nabla E\approx\lambda_{[x]}\,K^{T}_{[x]}\ast K_{[x]}\ast(u^{n}-g)-c_{[x]}\,\left(\nabla\cdot\nabla u^{n}-u_{\eta\eta}\right)+d_{[x]}\,u_{\eta\eta}\label{eq:gradient-local}
\end{equation}
where the subscript ${[x]}$ denotes the local point of spatially constant
approximation (rather than a function argument).
Assuming a uniform Cartesian grid oriented such that its first basis
vector $\vec{e_{1}}=(1,0,0,\ldots)$ aligns with $\nabla u$ at our
local point $x$, and that our spatial derivative discretizations
become equivalent to central difference (second derivative) approximations
with space step $\Delta x$ in each direction, then we obtain the
following local approximation of the gradient amplifier of \eqref{eq:gradient-local}
\begin{eqnarray}
  z(x,\omega_{1},\omega_{2},\ldots,\omega_{N})&\approx&
  \lambda_{[x]}\,\mbox{DFT}(K^{T}_{[x]})\,\mbox{DFT}(K_{[x]}) \nonumber \\
  &&+\frac{2}{\Delta x^{2}}\left(d_{[x]}\,(1-\cos\omega_{1}\Delta x)
  +c_{[x]}\sum_{k=2}^{N}(1-\cos\omega_{k}\Delta x)\right)
  \label{eq:amplifier-local}
\end{eqnarray}
Noting that the Fourier transform of the adjoint $K^{T}$ of a real
convolution kernel is always the complex conjugate of the Fourier
transform of the kernel $K$ itself, we see that the gradient amplifier
is real and positive and we can write the following upper bound as
a function frequency $\omega$
\[
\max_{\omega}z\le\lambda_{[x]}\max_{\omega}|\mbox{DFT}(K_{[x]})|^{2}+4\,\frac{c_{[x]}\,(N-1)+d_{[x]}}{\Delta x^{2}}
\]
with equality in cases where the complex amplitude of DFT($K$) is
maximal at $\omega=(\pi,\ldots,\pi)$. However, since this upper bound
depends on the local point of approximation $x$, we need to maximize
over $x$ as well in order to exploit the CFL formulas presented earlier
in terms of $z_{\max}$. Doing so yields the following upper bound
for the local gradient amplifier.
\begin{align}
z_{\max}\le & \max_{x,\omega}\left(|\mbox{DFT}(K)|^{2}\right)\lambda_{\max}+4\,\frac{(N-1)c_{\max}+d_{\max}}{\Delta x^{2}}\label{eq:zmax-general}\\
\mbox{where } & \lambda_{\max}\doteq\max_{x}\lambda,\quad c_{\max}\doteq\max_{x}c,\quad d_{\max}\doteq\max_{x}d\nonumber
\end{align}

If we now plug \eqref{eq:zmax-general} into the time step restriction
\eqref{eq:CFL-central} for the fully explicit second order accelerated
scheme \eqref{eq:FDE-central}, we obtain the following sufficient
condition for stability
\begin{equation}
\Delta t\le\frac{2\Delta x}{\sqrt{\max\left(|\mbox{DFT}(K)|^{2}\right)\lambda_{\max}\Delta x^{2}+4(N-1)c_{\max}+4d_{\max}}}\label{eq:CFL-general}
\end{equation}
The corresponding condition for gradient descent is obtained by squaring
$\Delta x$ in the numerator and removing the radical (squaring) the
denominator. As such we note three favorable step size trends for
PDE acceleration compared to PDE gradient descent. Most notably, when
the regularizing coefficients $c_{\max}$ and $d_{\max}$ dominate,
stable time step sizes are now directly proportional to spatial step
sizes rather than to their squares, making the upper bound linear
rather than quadratic in $\Delta x$. We see similar gains as well
when the kernel $K$ exhibits large amplification at one or more frequencies.
In such cases, stable step sizes are inversely proportional to the
maximum kernel amplification rather than to its square.

\subsection{Quadratic regularization (linear wave equation)\label{sec:quadratic}}

The easiest special case to consider would be that of quadratic fidelity
and regularity penalties without any forward model (more precisely,
with ${\cal K}$ as the identity operator).

\[
E(u)=\int_{\Omega}\frac{\lambda}{2}\left(u-g\right)^{2}+\frac{c}{2}\,\|\nabla u\|^{2}\,dx
\]
In this case the gradient is linear and the local approximation \eqref{eq:gradient-local}
becomes exact with $\lambda(x)=\lambda$, $c(x)=d(x)=c$.
\[
\nabla E=\lambda(u-g)-c\,\nabla\cdot\nabla u
\]
The accelerated descent PDE therefore takes the form of a damped inhomogeneous
linear wave equation.
\begin{equation}
u_{tt}-c\,\nabla\cdot\nabla u+au_{t}=\lambda(g-u)
\end{equation}
In this case the gradient amplifier $z(\omega)$ \eqref{eq:grad_amplifier}
is easy to compute. If central differences on a uniform $N$-dimensional
Cartesian grid with space step $\Delta x$ in each direction are used
to approximate the spatial derivatives of the Laplacian $\nabla\cdot\nabla$,
then
\[
z(\omega)=\lambda+\frac{2c}{\Delta x^{2}}\sum_{k=1}^{N}(1-\cos\omega_{k}\Delta x),\quad\omega=(\omega_{1},\ldots,\omega_{N})
\]
which makes the local approximation \eqref{eq:amplifier-local} exact
as well. Its upper bound
\begin{equation}
z_{\max}=\lambda+\frac{4Nc}{\Delta x^{2}}\label{eq:zmax-quadratic}
\end{equation}
is attained at $\omega=(\pi,\ldots,\pi)$, thereby making the general
condition \eqref{eq:CFL-general} necessary as well as sufficient
for stability. Plugging all this into \eqref{eq:update-gradient},
\eqref{eq:update-forward}, and \eqref{eq:update-central} yields
the following fully explicit updates (and CFL conditions), with multi-index
$\alpha=(\alpha_{1},\alpha_{2},\ldots,\alpha_{N})$ to indicate each
grid location, and where the additive multi-index $e_{k}=(\delta_{1k},\delta_{2k},\ldots,\delta_{Nk})$
is used to denote displacements to adjacent grid neighbors ($\delta_{jk}$
being the standard Kronecker delta).
\begin{align}
\mbox{\footnotesize gradient descent} & \begin{cases}
\Delta u_{\alpha}^{n}=-\Delta t\,\left(\lambda(u_{\alpha}^{n}-g_{\alpha})-c\,\sum_{k=1}^{N}\frac{u_{\alpha+e_{k}}^{n}-2u_{\alpha}^{n}+u_{\alpha-e_{k}}^{n}}{\Delta x^{2}}\right)\\
_{\Delta t\le\frac{2\Delta x^{2}}{4Nc+\lambda\Delta x^{2}}}
\end{cases}\label{eq:quadratic-gradient}\\
\mbox{\footnotesize 1-order accelerated} & \begin{cases}
\Delta u_{\alpha}^{n}=\frac{1}{1+a\Delta t}\Delta u_{\alpha}^{n-1}-\frac{\Delta t^{2}}{1+a\Delta t}\left(\lambda(u_{\alpha}^{n}-g_{\alpha})-c\,\sum_{k=1}^{N}\frac{u_{\alpha+e_{k}}^{n}-2u_{\alpha}^{n}+u_{\alpha-e_{k}}^{n}}{\Delta x^{2}}\right)\\
_{\Delta t\le\Delta x\sqrt{\frac{4}{4Nc+\lambda\Delta x^{2}}+\left(\frac{a\Delta x}{4Nc+\lambda\Delta x^{2}}\right)^{2}}+\frac{a\Delta x^{2}}{4Nc+\lambda\Delta x^{2}}}
\end{cases}\label{eq:quadratic-forward}\\
\mbox{\footnotesize 2-order accelerated} & \begin{cases}
\Delta u_{\alpha}^{n}=\frac{2-a\Delta t}{2+a\Delta t}\Delta u_{\alpha}^{n-1}-\frac{2\Delta t^{2}}{2+a\Delta t}\left(\lambda(u_{\alpha}^{n}-g_{\alpha})-c\,\sum_{k=1}^{N}\frac{u_{\alpha+e_{k}}^{n}-2u_{\alpha}^{n}+u_{\alpha-e_{k}}^{n}}{\Delta x^{2}}\right)\\
_{\Delta t\le\frac{2\Delta x}{\sqrt{4Nc+\lambda\Delta x^{2}}}}
\end{cases}\label{eq:quadratic-central}\\
\mbox{\footnotesize semi-implicit} & \begin{cases}
v_{\alpha}^{n}=u_{\alpha}^{n}+\frac{2-a\Delta t}{2+a\Delta t}\Delta u^{n-1}\\
u_{\alpha}^{n+1}=v_{\alpha}^{n}-\frac{2\Delta t^{2}}{2+a\Delta t}\,\left(\lambda(v_{\alpha}^{n}-g_{\alpha})-c\,\sum_{k=1}^{N}\frac{v_{\alpha+e_{k}}^{n}-2v_{\alpha}^{n}+v_{\alpha-e_{k}}^{n}}{\Delta x^{2}}\right)\\
_{\Delta t\le\frac{2\Delta x}{\sqrt{3\left(4Nc+\lambda\Delta x^{2}\right)}}\quad(\mbox{sufficent but not necessary when \ensuremath{a>0})}}
\end{cases}\label{eq:quadratic-semi}
\end{align}

\subsection{Implicit handling of the fidelity term\label{sec:implicit}}

The portion of the continuum gradient which arises from the fidelity
term is $\lambda(u-g)$, which we have discretized explicitly in the
above schemes as $\lambda(u_{\alpha}^{n}-g_{\alpha})$. Since this
term, unlike the Laplacian discretization, does not depend upon neighboring
grid locations, we could evaluate it implicitly at the updated value
of $u$ by plugging $\lambda(u_{\alpha}^{n+1}-g_{\alpha})$ into any
of these schemes and yet still rearrange the resulting expressions
to obtain explicit updates for $u_{\alpha}^{n+1}$. Algebraic manipulation
of these resulting implicitly handled fidelity schemes would yield
the following equivalent schemes, restructured to reveal their similarity
to the schemes \eqref{eq:quadratic-gradient}, \eqref{eq:quadratic-forward},
\eqref{eq:quadratic-central}, and \eqref{eq:quadratic-semi} shown
above.
\begin{align}
\mbox{\footnotesize gradient descent: } & \Delta u_{\alpha}^{n}=-\frac{\Delta t}{1+\lambda\Delta t}\left(\lambda(u_{\alpha}^{n}-g_{\alpha})-c\,\sum_{k=1}^{N}\frac{u_{\alpha+e_{k}}^{n}-2u_{\alpha}^{n}+u_{\alpha-e_{k}}^{n}}{\Delta x^{2}}\right)\label{eq:implicit-gradient}\\
\mbox{\footnotesize 1-order accelerated: } & \Delta u_{\alpha}^{n}=\frac{\Delta u_{\alpha}^{n-1}-\Delta t^{2}\left(\lambda(u_{\alpha}^{n}-g_{\alpha})-c\,\sum_{k=1}^{N}\frac{u_{\alpha+e_{k}}^{n}-2u_{\alpha}^{n}+u_{\alpha-e_{k}}^{n}}{\Delta x^{2}}\right)}{1+\left(a+\lambda\Delta t\right)\Delta t}\label{eq:implicit-forward}\\
\mbox{\footnotesize 2-order accelerated: } & u_{\alpha}^{n}=\frac{\left(2-a\Delta t\right)\Delta u_{\alpha}^{n-1}-2\Delta t^{2}\left(\lambda(u_{\alpha}^{n}-g_{\alpha})-c\,\sum_{k=1}^{N}\frac{u_{\alpha+e_{k}}^{n}-2u_{\alpha}^{n}+u_{\alpha-e_{k}}^{n}}{\Delta x^{2}}\right)}{2+a\Delta t+2\lambda\Delta t^{2}}\label{eq:implicit-central}\\
\mbox{\footnotesize semi-implicit: } & \begin{cases}
v_{\alpha}^{n}=u_{\alpha}^{n}+\frac{2-a\Delta t}{2+a\Delta t}\Delta u^{n-1}\\
u_{\alpha}^{n+1}=v_{\alpha}^{n}-\frac{2\Delta t^{2}}{2+a\Delta t+2\lambda\Delta t^{2}}\,\left(\lambda(v_{\alpha}^{n}-g_{\alpha})-c\,\sum_{k=1}^{N}\frac{v_{\alpha+e_{k}}^{n}-2v_{\alpha}^{n}+v_{\alpha-e_{k}}^{n}}{\Delta x^{2}}\right)
\end{cases}\label{eq:implicit-semi}
\end{align}
Written in this form it is easy to show by comparison that these schemes
become equivalent to their explicit-fidelity counterparts by a change
of time step, damping parameter, or both. In the case of gradient
descent, the implicit-fidelity scheme \eqref{eq:implicit-gradient}
is identical to explicit-fidelity scheme \eqref{eq:quadratic-gradient}
with a smaller time step, using $\Delta t\to\frac{\Delta t}{1+\lambda\Delta t}$.
The first order implicit-fidelity accelerated scheme \eqref{eq:implicit-forward}
is equivalent to its explicit-fidelity counterpart \eqref{eq:quadratic-forward}
with a larger damping coefficient, using $a\to a+\lambda\Delta t$.
The second order implicit-fidelity accelerated scheme \eqref{eq:implicit-central}
is equivalent to the explicit-fidelity scheme \eqref{eq:quadratic-central}
with both a smaller time step and an adjusted damping coefficient
(may be either larger or smaller depending on $\lambda$), using $\Delta t\to\frac{\Delta t}{\sqrt{1+\frac{\lambda}{2}\Delta t^{2}}}$
and $a\to\frac{a+\lambda\Delta t}{\sqrt{1+\frac{\lambda}{2}\Delta t^{2}}}$.
Finally the implicit-fidelity adaptation \eqref{eq:implicit-semi}
of the semi-implicit scheme \eqref{eq:quadratic-semi}, obtained by
replacing $\lambda(u_{\alpha}^{n+1}-g)$ with $\lambda(v_{\alpha}^{n}-g)$,
is equivalent to the original semi-implicit scheme \eqref{eq:quadratic-semi}
with both a smaller time step and a larger damping coefficient, using
$\Delta t\to\frac{\Delta t}{\sqrt{1+\frac{2\lambda\Delta t^{2}}{2+a\Delta t}}}$
and $a\to a\sqrt{1+\frac{2\lambda\Delta t^{2}}{2+a\Delta t}}$.

The CFL conditions for these implicit-fidelity schemes can therefore
be obtained by applying these substitutions backwards to the matching
explicit (or semi-implicit) CFL conditions. While this often yields
a larger maximum stable time step, the apparent gain is deceptive
since there is will be no numerical difference to the corresponding
explicit (or semi-implicit) update with a smaller time step. As such,
there is neither a computational nor a numerical advantage to handling
the fidelity term implicitly. While we have illustrated this here
for the special case of quadratic regularization, the parameter remappings
showing equivalency between the explicit and partially implicit schemes
depend only upon the damping and fidelity parameters. It is easy to
see that the exact same analysis applies even in the nonlinear case
of non-quadratic regularization, making this equivalency (and therefore
the lack of benefit in implicitly handling the fidelity) more general.

Further generalization of this analysis is also possible in the accelerated
cases for non-quadratic fidelity penalization as well as for nontrivial
forward models ${\cal K}$. However, in such cases, equivalency would
require substitution of a constant damping parameter $\alpha$ in
the partially implicit scheme with a spatially varying damping in
the equivalent explicit scheme. For example, in the case of a quadratic
fidelity penalty paired with a convolution kernel $K$ in the first
order accelerated implicit-fidelity scheme \eqref{eq:implicit-forward},
a constant damping parameter $a$ would be have to be replaced by
the spatially varying $a{\cal I}+\lambda\Delta tK^{^{T}}K$ in order
to use the explicit-fidelity scheme \eqref{eq:quadratic-forward}
to obtain equivalent updates. This would require inversion of the
matrix $(1+a\Delta t){\cal I}+\lambda\Delta t^{2}K^{^{T}}K$, as division
by a scalar would no longer occur in the explicit update \eqref{eq:quadratic-forward}.
However since this inverse does not depend on $u$, it inverse could
be computed/approximated just once and then reused in every update
step (in cases where the damping does not change with time).

\subsection{Beltrami regularization (quasi-linear wave equation)
\label{sec:beltrami}}

Another special case to consider is Beltrami regularization.
We'll consider the case of a quadratic penalty and an attenuating,
mean-preserving convolution kernel $K$

\begin{equation} \label{eq:beltrami-energy}
E(u)=\int_{\Omega}\frac{\lambda}{2}\left(K\ast u-g\right)^{2}+\underbrace{\frac{1}{\beta}\sqrt{1+\|\beta\nabla u\|^{2}}}_{\sqrt{\epsilon^{2}+\|\nabla u\|^{2}},\;\epsilon=\frac{1}{\beta}}\,dx
\end{equation}
In this case the variational gradient is non-linear and \eqref{eq:gradient-local}
decomposes as follows.
\begin{align*}
\nabla E & =\lambda K^{T}\ast K\ast(u-g)-\nabla\cdot\underbrace{\left(\frac{\beta\nabla u}{\sqrt{1+\|\beta\nabla u\|^{2}}}\right)}_{\frac{\nabla u}{\sqrt{\epsilon^{2}+\|\nabla u\|^{2}}},\;\epsilon=\frac{1}{\beta}}\\
 & =\lambda K^{T}\ast K\ast(u-g)-\underbrace{\frac{\beta}{\sqrt{1+\|\beta\nabla u\|^{2}}}}_{c}\left(\nabla\cdot\nabla u-u_{\eta\eta}\right)\,-\,\underbrace{\frac{\beta}{\left(\sqrt{1+\|\beta\nabla u\|^{2}}\right)^{3}}}_{d}\;u_{\eta\eta}
\end{align*}
while the accelerated PDE (technically an integral partial differential
equation with the convolution) takes the quasilinear form.

\begin{equation} \label{eq:beltrami-PDE}
u_{tt}-\nabla\cdot\left(\frac{\beta\nabla u}{\sqrt{1+\|\beta\nabla u\|^{2}}}\right)+au_{t}=\lambda K^{T}\ast K\ast(g-u)
\end{equation}
Note that both coefficients $c$ and $d$ are bounded by $\beta$
(an upper bound which is actually reached in both cases at any point
and time where $\nabla u(x,t)=0$), and that $\max\left|\mbox{DFT}(K)\right|=1$
by our assumption that $K$ attenuates while preserving the mean.
Plugging this into \eqref{eq:zmax-general} yields
\begin{equation}
z_{\max}\le\lambda+\frac{4N\beta}{\Delta x^{2}}\label{eq:zmax-beltrami}
\end{equation}
if we assume a consistent discretization of $\nabla\cdot\left(\frac{\beta\nabla u}{\sqrt{1+\|\beta\nabla u\|^{2}}}\right)$
which converges, as $\nabla u\to0$, to the central difference approximation
of the $\beta$-scaled Laplacian $\beta\nabla\cdot\nabla u\approx\beta\sum_{k=1}^{N}\frac{u_{\alpha+e_{k}}^{n}-2u_{\alpha}^{n}+u_{\alpha-e_{k}}^{n}}{\Delta x^{2}}$
with spatial step size $\Delta x$ in each direction (see Section
\eqref{sec:quadratic} for the multi-index subscript notation $\alpha$
and $e_{k}$). If we let $D_{\!\beta,\Delta x}^{2}u^{n}$ denote the
discretization of $\nabla\cdot\left(\frac{\beta\nabla u}{\sqrt{1+\|\beta\nabla u\|^{2}}}\right)$
then we obtain the following schemes
\begin{align}
\mbox{\footnotesize gradient descent} & \begin{cases}
\Delta u^{n}=-\Delta t\,\left(\lambda K^{T}\ast K\ast(u^{n}-g)-D_{\!\beta,\Delta x}^{2}u^{n}\right)\\
_{\Delta t\le\Delta x^{2}\left(\frac{2}{4N\beta+\lambda\Delta x^{2}}\right)}
\end{cases}\label{eq:beltrami-gradient}\\
\mbox{\footnotesize 1-order accelerated} & \begin{cases}
\Delta u^{n}=\frac{1}{1+a\Delta t}\Delta u^{n-1}-\frac{\Delta t^{2}}{1+a\Delta t}\left(\lambda K^{T}\ast K\ast(u^{n}-g)-D_{\!\beta,\Delta x}^{2}u^{n}\right)\\
_{\Delta t\le\Delta x\left(\sqrt{\frac{4}{4N\beta+\lambda\Delta x^{2}}+\left(\frac{a\Delta x}{4N\beta+\lambda\Delta x^{2}}\right)^{2}}+\frac{a\Delta x}{4N\beta+\lambda\Delta x^{2}}\right)}
\end{cases}\label{eq:beltrami-forward}\\
\mbox{\footnotesize 2-order accelerated} & \begin{cases}
\Delta u^{n}=\frac{2-a\Delta t}{2+a\Delta t}\Delta u^{n-1}-\frac{2\Delta t^{2}}{2+a\Delta t}\left(\lambda K^{T}\ast K\ast(u^{n}-g)-D_{\!\beta,\Delta x}^{2}u^{n}\right)\\
_{\Delta t\le\Delta x\left(\frac{2}{\sqrt{4N\beta+\lambda\Delta x^{2}}}\right)}
\end{cases}\label{eq:beltrami-central}\\
\mbox{\footnotesize semi-implicit} & \begin{cases}
v^{n}=u^{n}+\frac{2-a\Delta t}{2+a\Delta t}\Delta u^{n-1}\\
u^{n+1}=v^{n}-\frac{2\Delta t^{2}}{2+a\Delta t}\,\left(\lambda K^{T}\ast K\ast(v^{n}-g)-D_{\!\beta,\Delta x}^{2}v^{n}\right)\\
_{\Delta t\le\Delta x\left(\frac{2}{\sqrt{3\left(4N\beta+\lambda\Delta x^{2}\right)}}\right)}
\end{cases}\label{eq:beltrami-semi}
\end{align}

\subsection{Total Variation Regularization\label{sec:TV}}

If we consider the limit as $\beta\to\infty$, the Beltrami regularization
penalty converges to the total variation penalty.

\begin{equation} \label{eq:TV-energy}
E(u)=\int_{\Omega}\frac{\lambda}{2}\left(K\ast u-g\right)^{2}+\|\nabla u\|\,dx
\end{equation}
with a non-linear variational gradient \eqref{eq:gradient-local}
that decomposes as follows.
\begin{align*}
\nabla E & =\lambda K^{T}\ast K\ast(u-g)-\nabla\cdot\left(\frac{\nabla u}{\|\nabla u\|}\right)=\lambda K^{T}\ast K\ast(u-g)-\underbrace{\frac{1}{\|\nabla u\|}}_{c}\left(\nabla\cdot\nabla u-u_{\eta\eta}\right)
\end{align*}
The accelerated PDE now takes the form of the following nonlinear
wave equation.
\begin{equation} \label{eq:TV-PDE}
u_{tt}-\nabla\cdot\left(\frac{\nabla u}{\|\nabla u\|}\right)+au_{t}=\lambda K^{T}\ast K\ast(g-u)
\end{equation}
In this case, the coefficient $d$ vanishes, but the coefficient $c$
no longer has a finite upper bound. Plugging this into \eqref{eq:zmax-general}
yields an infinite upper bound for the maximum gradient amplifier
if at any point and time $\nabla u(x,t)=0$. Otherwise, by our earlier
assumption on $K$ (see Section \ref{sec:beltrami}) we obtain
\begin{equation}
\lambda\le
z_{\max}
\le\lambda+\frac{4(N-1)}{\Delta x^{2}\min\|\nabla u\|}.\label{eq:zmax-TV}
\end{equation}

For the explicit second order accelerated scheme, this ensures the
sufficient condition $\Delta t\le\frac{2}{\sqrt{\lambda+\frac{4(N-1)}{\Delta x^{2}\min\|\nabla u\|}}}$
for a stable step. If we fix $\Delta t$, we may rearrange this inequality
to obtain an equivalent sufficient condition
\[
\min\|\nabla u\|\ge\frac{N-1}{\Delta x^{2}}\frac{4\Delta t^{2}}{4-\lambda\Delta t^{2}}
\]
which takes the form of a lower bound on the spatial gradient.

Here
an interesting nonlinear dynamic occurs to keep the implementation
stable by preventing initiated instabilities from growing unbounded.
If the spatial gradient falls below this lower bound and instabilities
begins to propagate at one or more frequencies, they will eventually
cause the spatial gradient to rise above the guaranteed stable lower
bound at which point the instabilities will cease growing. In the
absence of a kernel $K$, the fastest growing instability will occur
at the highest digital frequency in each grid direction $\omega=(\pi,\ldots,\pi)$
which corresponds to oscillations between immediately adjacent grid-points,
this in turn will most rapidly increase the discrete difference approximations
of $\|\nabla u\|$. In the presence of a strongly smoothing kernel,
the fastest growing instability may occur at lower digital frequencies,
thereby causing a low-grade ringing effect, with several grid-points
per period, until the amplitude of the oscillation is large enough
to drive adjacent pixel differences back over the lower bound for
$\|\nabla u\|$.

A similar phenomenon occurs with both the first-order and
semi-implicit schemes (and even with gradient descent), making all these
schemes stable independently of the regularizer coefficient $c$. As such,
purely for stability considerations alone, the necessary step size constraint
will be connected to the lower bound $\lambda$ of the gradient amplifier
$z_{\max}$ rather than its upper bound in \eqref{eq:zmax-TV}. This yields
the following necessary conditions for stability.
\begin{align}
\mbox{\footnotesize gradient descent: }	  &\label{eq:TV-gradient-CFL}\Delta t\le
\frac{2}{\lambda} \\
\mbox{\footnotesize 1-order accelerated: }&\label{eq:TV-forward-CFL}\Delta t\le
\sqrt{\frac{4}{\lambda}+\left(\frac{a}{\lambda}\right)^{2}}+\frac{a}{\lambda} \\
\mbox{\footnotesize 2-order accelerated: }&\label{eq:TV-central-CFL}\Delta t\le
\frac{2}{\sqrt{\lambda}} \\
\mbox{\footnotesize semi-implicit: }	  &\label{eq:TV-semi-CFL}\Delta t\le
\frac{2}{\sqrt{3\lambda}}
\end{align}
However, the schemes may only converge under these constraints in an
oscillatory sense with a fluctuating level of ``background noise'' whose
amplitude will depend upon the value of $\Delta t$.

We may exploit the behavior of this non-linear stabilizing effect
to obtain a more useful time step constraint by plugging in a
minimal acceptable value of $\|\nabla u\|$ for the final reconstruction
into the stability condition for $\Delta t$. A natural way to approach
this is by exploiting a quantization interval $Q$ for the digital
representation of $u$ together with the following discrete approximation
bounds for $||\nabla u\|$.
\[
\min\|\nabla u\|=\min_{\alpha}\sqrt{\sum_{k=1}^{N}\left(\frac{u_{\alpha+e_{k}}-u_{\alpha}}{\Delta x}\right)^{2}}\ge\sqrt{N\min_{\alpha,k}\left(\frac{u_{\alpha+e_{k}}-u_{\alpha}}{\Delta x}\right)^{2}}=\frac{\sqrt{N}}{\Delta x}\min_{\alpha,k}\left|u_{\alpha+e_{k}}-u_{\alpha}\right|
\]
If we now determine that instability related distortions confined
to a single quantization interval $Q$ between neighboring pixels
are acceptable, we substitute
\[
\min\|\nabla u\|\to\frac{\sqrt{N}}{\Delta x}Q
\]
 into the upper bound for \eqref{eq:zmax-TV} to obtain
\begin{equation}
z_{\max}\le\lambda+\frac{4(N-1)}{Q\Delta x\sqrt{N}}<\lambda+\frac{4\sqrt{N}}{Q\Delta x}\label{eq:zmax-TV-quantized}
\end{equation}
within the desired stable regime for $\|\nabla u\|$. This in turn
gives rise to the following schemes, where $D_{\!\Delta x}^{2}u^{n}$
denotes the discretization of $\nabla\cdot\left(\frac{\nabla u}{\|\nabla u\|}\right)$
\begin{align}
\mbox{\footnotesize gradient descent} & \begin{cases}
\Delta u^{n}=-\Delta t\,\left(\lambda K^{T}\ast K\ast(u^{n}-g)-D_{\!\Delta x}^{2}u^{n}\right)\\
_{\Delta t\le Q\Delta x\,\left(\frac{2}{4\sqrt{N}+\lambda Q\Delta x}\right)}
\end{cases}\label{eq:TV-gradient}\\
\mbox{\footnotesize 1-order accelerated} & \begin{cases}
\Delta u^{n}=\frac{1}{1+a\Delta t}\Delta u^{n-1}-\frac{\Delta t^{2}}{1+a\Delta t}\left(\lambda K^{T}\ast K\ast(u^{n}-g)-D_{\!\Delta x}^{2}u^{n}\right)\\
_{\Delta t\le\sqrt{Q\Delta x\,}\left(\sqrt{\frac{4}{4\sqrt{N}+\lambda Q\Delta x}+\left(\frac{a\sqrt{Q\Delta x}}{4\sqrt{N}+\lambda Q\Delta x}\right)^{2}}+\frac{a\sqrt{Q\Delta x}}{4\sqrt{N}+\lambda Q\Delta x}\right)}
\end{cases}\label{eq:TV-forward}\\
\mbox{\footnotesize 2-order accelerated} & \begin{cases}
\Delta u^{n}=\frac{2-a\Delta t}{2+a\Delta t}\Delta u^{n-1}-\frac{2\Delta t^{2}}{2+a\Delta t}\left(\lambda K^{T}\ast K\ast(u^{n}-g)-D_{\!\Delta x}^{2}u^{n}\right)\\
_{\Delta t\le\sqrt{Q\Delta x\,}\left(\sqrt{\frac{4}{4\sqrt{N}+\lambda Q\Delta x}}\right)}
\end{cases}\label{eq:TV-central}\\
\mbox{\footnotesize semi-implicit} & \begin{cases}
v^{n}=u^{n}+\frac{2-a\Delta t}{2+a\Delta t}\Delta u^{n-1}\\
u^{n+1}=v^{n}-\frac{2\Delta t^{2}}{2+a\Delta t}\,\left(\lambda K^{T}\ast K\ast(v^{n}-g)-D_{\!\Delta x}^{2}v^{n}\right)\\
_{\Delta t\le\sqrt{Q\Delta x\,}\left(\sqrt{\frac{4}{3\left(4\sqrt{N}+\lambda Q\Delta x\right)}}\right)}
\end{cases}\label{eq:TV-semi}
\end{align}

\section{Experimental Examples} \label{sec:experiments}

\subsection{Beltrami Denoising}

\begin{figure}
\centering
\begin{tabular}{c@{}c@{}c@{}c@{}c}
\epsfig{bb=120 90 274 448,clip=true,width=0.20\textwidth,figure=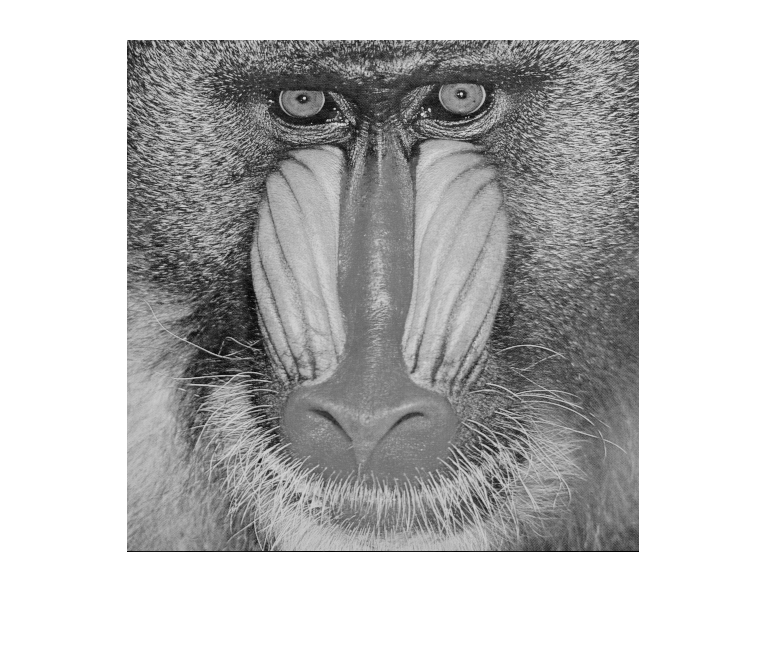} &
\epsfig{bb=120 90 274 448,clip=true,width=0.20\textwidth,figure=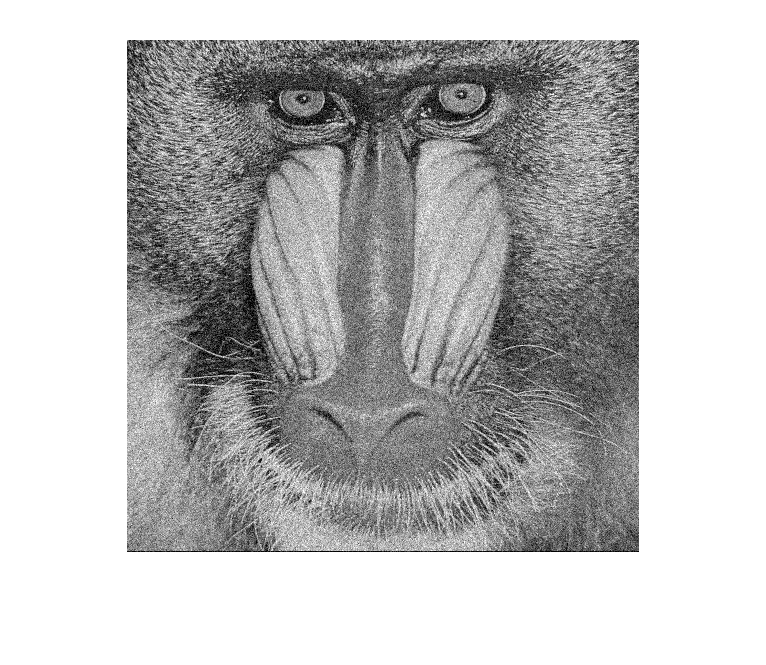} &
\epsfig{bb=120 90 274 448,clip=true,width=0.20\textwidth,figure=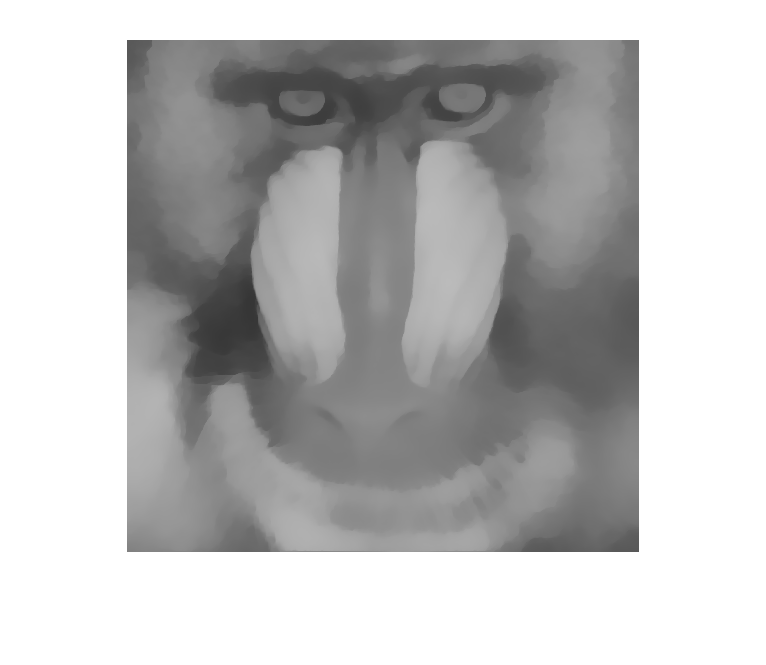} &
\epsfig{bb=120 90 274 448,clip=true,width=0.20\textwidth,figure=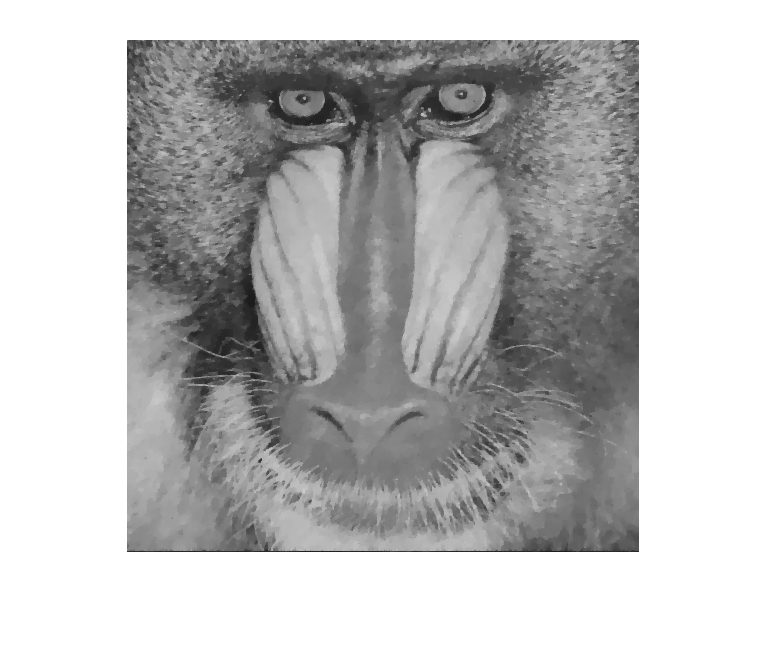} &
\epsfig{bb=120 90 274 448,clip=true,width=0.20\textwidth,figure=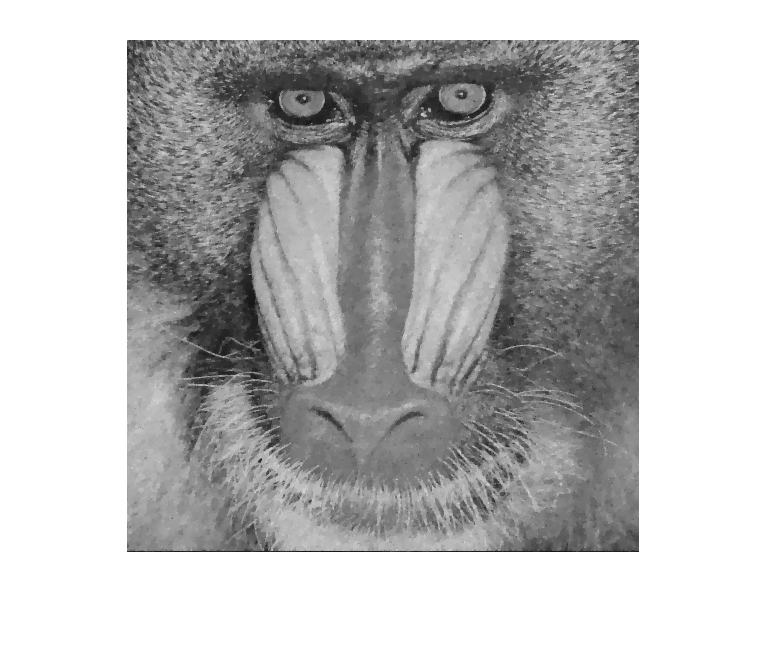} \\
Original & Noisy &$\beta^2=1,\lambda=1$ & $\beta^2=1,\lambda=5$ & $\beta^2=1,\lambda=7$ \\
\epsfig{bb=120 90 274 448,clip=true,width=0.20\textwidth,figure=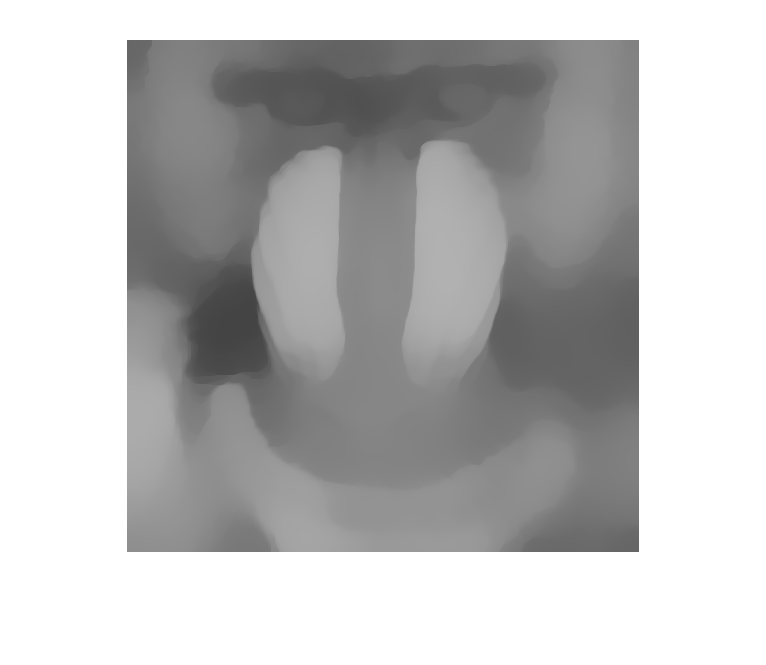} &
\epsfig{bb=120 90 274 448,clip=true,width=0.20\textwidth,figure=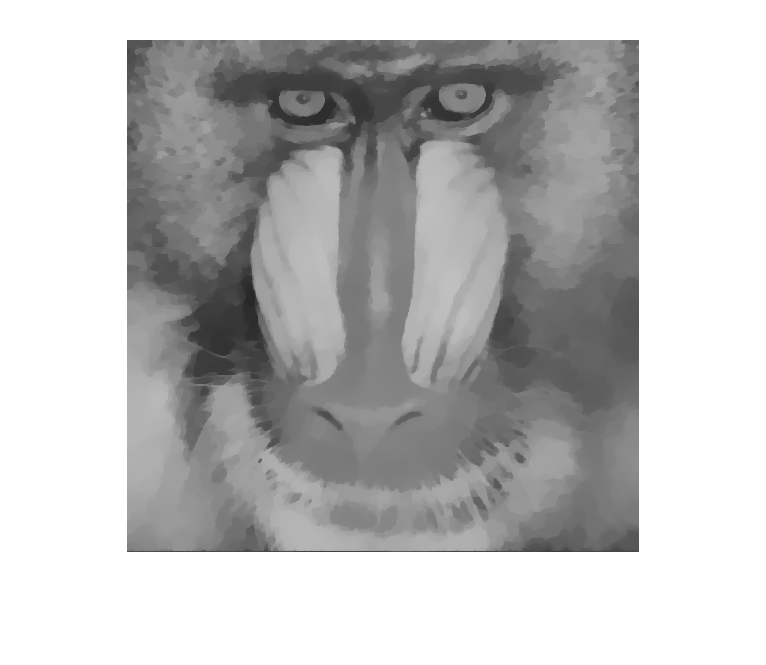} &
\epsfig{bb=120 90 274 448,clip=true,width=0.20\textwidth,figure=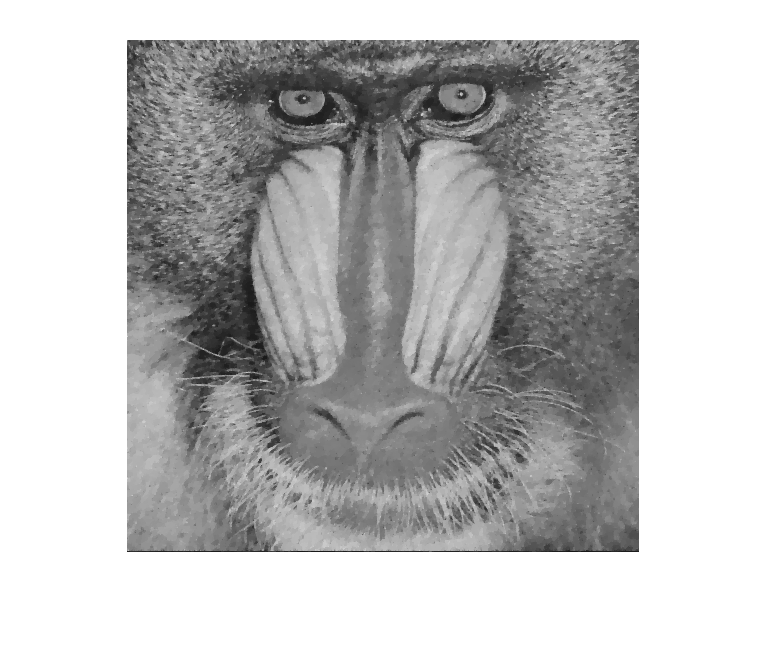} &
\epsfig{bb=120 90 274 448,clip=true,width=0.20\textwidth,figure=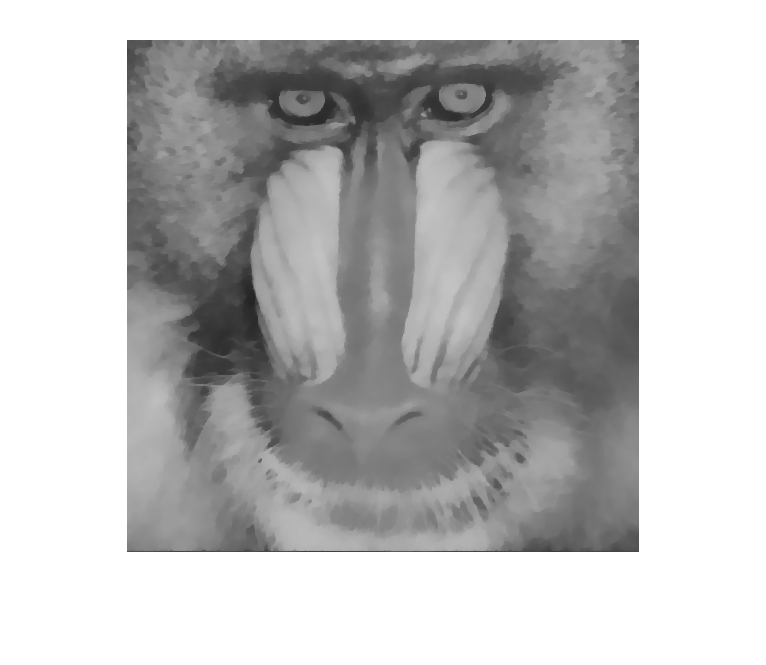} &
\epsfig{bb=120 90 274 448,clip=true,width=0.20\textwidth,figure=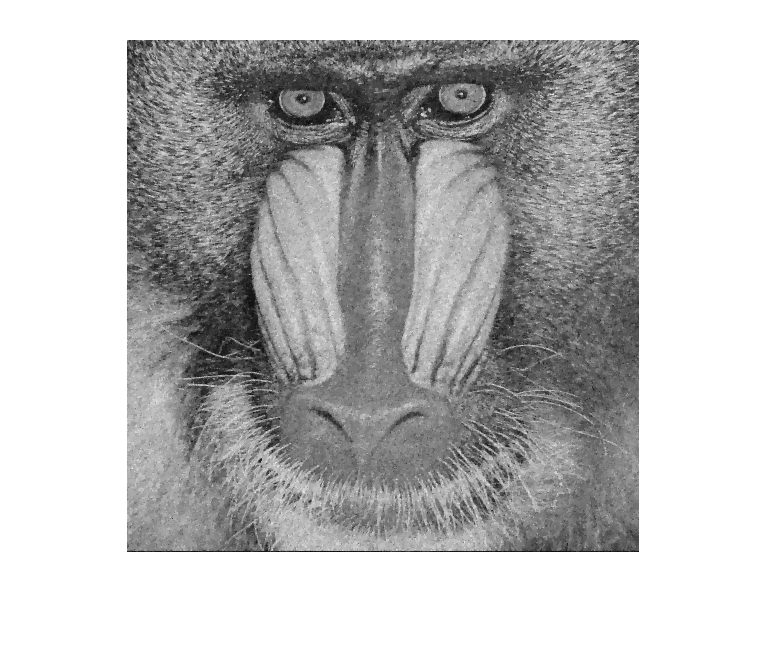} \\
$\beta^2=5,\lambda=0.5$ & $\beta^2=5,\lambda=2$ & $\beta^2=5,\lambda=7$ & $\beta^2=1/5,\lambda=2$ & $\beta^2=1/5,\lambda=10$ \\
\end{tabular}
\caption{Results of Beltrami regularization applied to a noisy baboon image with varying values of $\lambda$ and $\beta$. The units of $\lambda$ are thousands.}
\label{fig:baboon}
\end{figure}

Our first application is to the problem of Beltrami regularization
for image denoising and image inpainting
\cite{zosso2014primal,kimmel1998image,sochen1998general},
which corresponds to minimizing \eqref{eq:beltrami-energy} in the absence
of a kernel $K$ via the accelerated PDE \eqref{eq:beltrami-PDE}. In this case
$g$ is the original noisy image,
and the minimizer $u$ is the denoised/inpainted
image. For denoising we typically set the parameter $\lambda$ to be a positive
constant, and for inpainting we can set $\lambda=0$ in the region
$D\subset\Omega$ to be inpainted, and set $\lambda$ to be large
or $\infty$ in $\Omega\setminus D$. The Beltrami regularization
term interpolates between the TV norm $\int|\nabla u|$ and the $H^{1}$
norm $\int|\nabla u|^{2}$---near edges where $\nabla u$ is large,
it behaves like the TV norm to preserve edges, and where $\nabla u$
is small it behaves like the $H^{1}$ norm in order to reduce staircasing.
Recently, Zosso and Bustin \cite{zosso2014primal} proposed an efficient
primal dual projected gradient method for solving Beltrami regularized
problems.

\begin{table}
\centering
\begin{tabular}{|c|c|c|c|c|c|c|}
 \hline
 &\multicolumn{2}{c|}{$\lambda=1000$}&\multicolumn{2}{c|}{$\lambda=5000$} &\multicolumn{2}{c|}{$\lambda=7000$}\\
\hline
 & Time & Iterations& Time & Iterations & Time & Iterations \\
\hline\hspace{2mm}
$\beta^2=1/5$   & 0.55s & 124       & 0.27s &  60       & 0.23s &  50 \\
$\beta^2=1$     & 0.81s & 183       & 0.38s &  85       & 0.32s &  71 \\
$\beta^2=5$     & 1.20s & 273       & 0.54s & 122       & 0.45s & 101 \\
\hline
\end{tabular}
\caption{PDE accelerated Beltrami regularization runtimes on the $512\times 512$ baboon image.}
\label{tab:BRsim}
\end{table}

We use the first order explicit scheme \eqref{eq:beltrami-forward}
with forward differences for $\nabla u$ and backward
differences for $\text{div}$.
We set the damping coefficient to $a=2\sqrt{\beta\pi^2+\lambda}$,
via the linear analysis \eqref{eq:aopt} and
run the algorithm at its maximum stable step $\Delta t$
\eqref{eq:beltrami-forward}
until the absolute
difference between the current and previous iterate
falls below $10^{-4}$. We note that the image is normalized so the
pixel values fall in the interval $[0,1]$.

Figure \ref{fig:baboon} shows the results of applying the PDE accelerated
Beltrami regularization to a noisy baboon test image with varying
values of $\lambda$ and $\beta$ with single-threaded C++ code on a
3.2 GHz Intel processor running Linux. The corresponding runtimes
are given in Table \ref{tab:BRsim} and are favorably competitive with
the runtimes reported in \cite{zosso2014primal},
who proposed a primal dual projected
gradient algorithm for Beltrami regularization. Notice the algorithm
does slow down somewhat when $\lambda$ is small and the denoising
is heavily regularized, but the difference is far less pronounced
compared to other explicit methods such as gradient descent.

\subsection{Beltrami Inpainting}

\begin{figure}
\centering
\subfloat[Original]{\includegraphics[trim = 90 90 100 50, clip=true,width=0.50\textwidth]{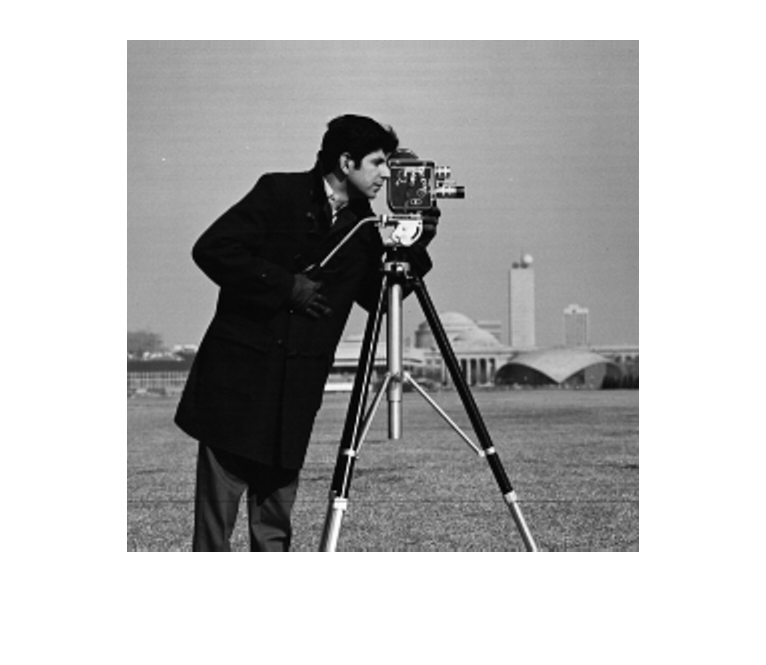}}
\subfloat[Corrupted]{\includegraphics[trim = 90 90 100 50, clip=true,width=0.50\textwidth]{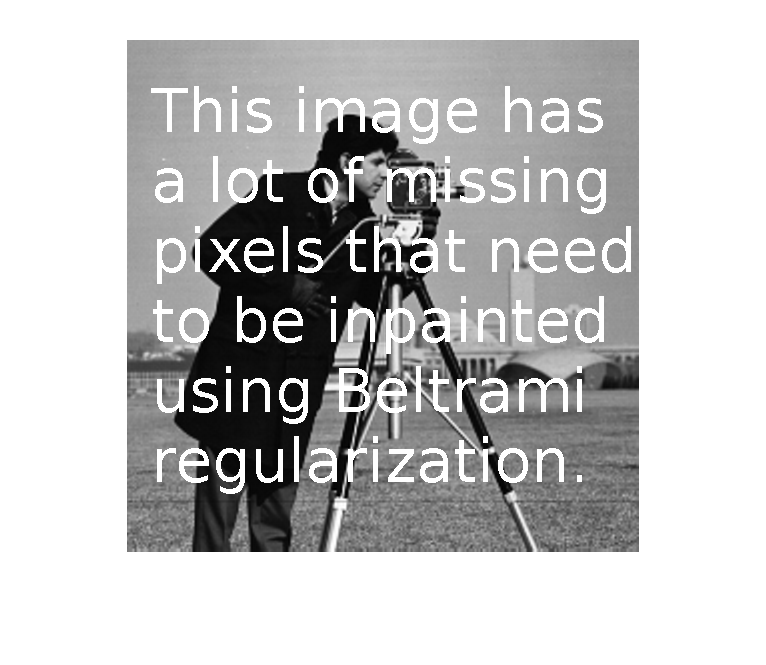}}

\subfloat[Inpainting result ($\beta=1$)]{\includegraphics[trim = 90 90 100 50, clip=true,width=0.50\textwidth]{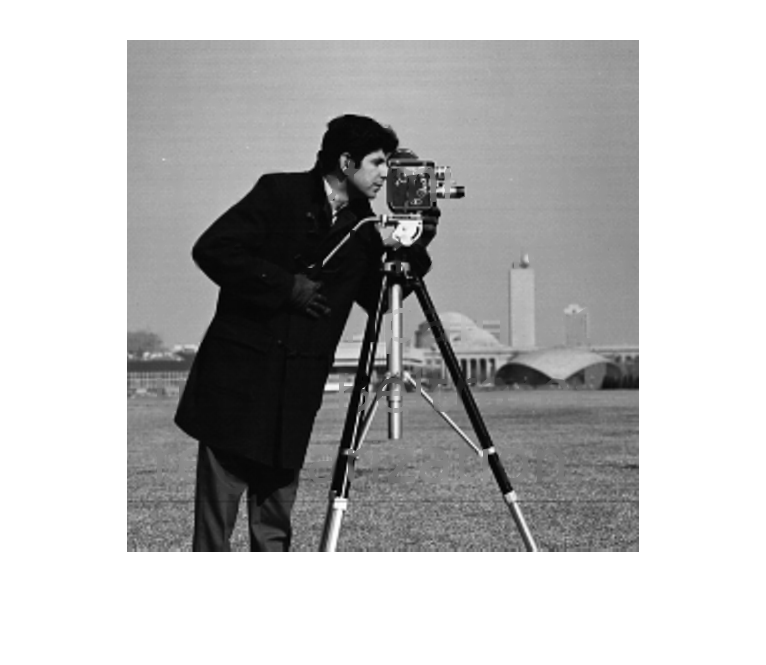}}
\hfill
\subfloat[Inpainting Enlarged]{\includegraphics[trim = 250 168.5 250 260, clip=true,width=0.495\textwidth]{cameramanIP}}
\caption{An example of inpainting using the PDE accelerated Beltrami regularization framework on the cameraman image.}
\label{fig:cameramanIP}
\end{figure}

We also give an example of PDE acceleration for Beltrami regularized
inpainting in Figure \ref{fig:cameramanIP}. We used $\beta=1$
and $a=5\pi$, and the inpainting took 687 iterations (11.48 seconds)
starting from an initial guess given by nearest neighbor interpolation.
This is a good deal slower than the denoising examples. It is possible
to give a partial explanation for this. Recall that the optimal damping
parameter, and convergence rate, depends on the size of the first
eigenvalue of the linearized operator on the given domain, and the
presence of a zeroth order term $\lambda u$. In inpainting, there
is no zeroth order term and the domain is highly irregular. Further,
the inpainting domain is typically disconnected, so the eigenvalues
on each connected component would be required, and this would lead
to different choices of damping coefficient in each region. We plan
to investigate this issue, and others, in future work.

\subsection{Beltrami Deblurring}

\begin{figure}
\centering
\begin{tabular}{c@{}c@{}c@{}c}
\epsfig{bb=0 10 500 510,clip=true,width=0.25\textwidth,figure=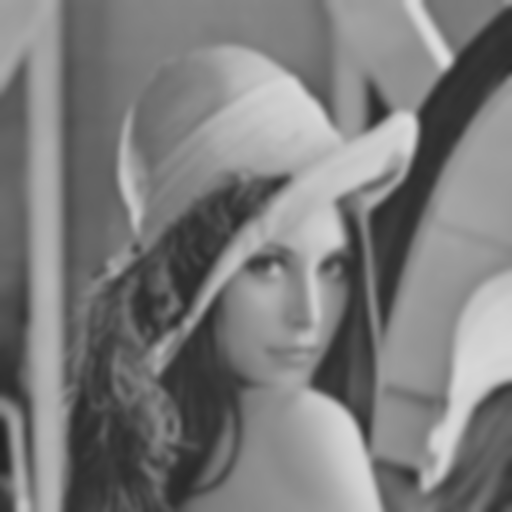} &
\epsfig{bb=0 10 500 510,clip=true,width=0.25\textwidth,figure=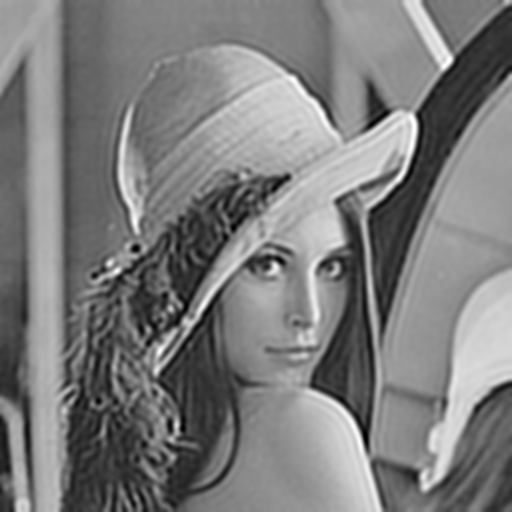} &
\epsfig{bb=0 10 500 510,clip=true,width=0.25\textwidth,figure=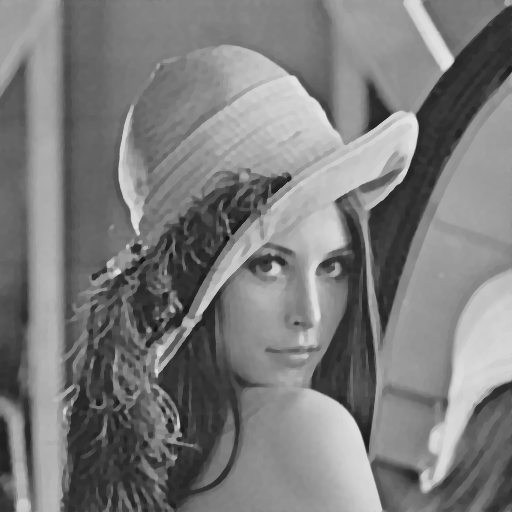} &
\epsfig{bb=0 10 500 510,clip=true,width=0.25\textwidth,figure=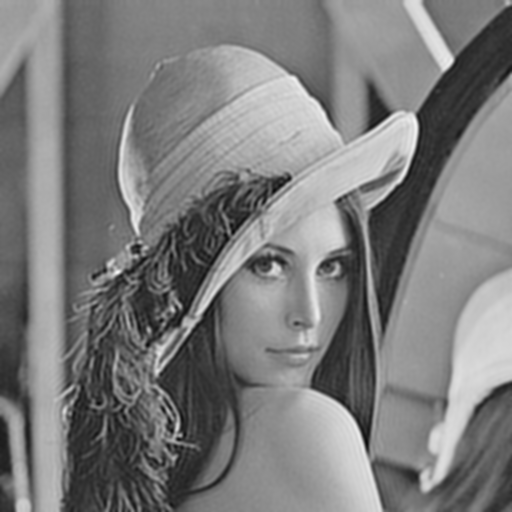} \\
Original & Primal Dual & L1 ADMM & PDE acceleration \vspace{-1mm}\\
\footnotesize PSNR=25.6dB &\footnotesize 27.8dB &\footnotesize 31.8dB &\footnotesize 32.3dB \vspace{-3mm}
\end{tabular}
\caption{\label{fig:leena_TVdeblur_fair}
Deblurring of an image using the explicit accelerated PDE scheme compared
with the results of two other state of the art methods (final
signal-to-noise ratios shown for each restoration).
}
\end{figure}

Finally, we give an example of PDE acceleration for Beltrami regularized
deblurring. We used $\lambda=10^7$,
$\beta=1$, and $a=4$ and the deblurring was run
using the second order explicit scheme \eqref{eq:beltrami-central} with its
maximum stable time step starting with the original blurred image as the
initial guess. After 2038 iterations it achieved its tenth-of-a-decibel 
rounded steady-state restored PSNR of 32.3dB.
The original image was blurred with a Gaussian kernel of $\sigma=3$ to create
an blurry initial image with a signal-to-noise ratio of 25.6185dB.
In Figure \ref{fig:leena_TVdeblur_fair} we compare the accelerated PDE results, both visually as well
as quantitatively according to the restored signal-to-noise ratio, with
those obtained using primal dual and L1 ADMM algorithms for the same parameters $\lambda=10^7$ and $\beta=1$. ADMM reached its tenth-of-a-decibel rounded
steady-state restored PSNR of 31.8dB after 2453 iterations, whereas Primal Dual
reached its tenth-of-a-decibel rounded steady-state restored PSNR of 27.8dB
after 63 iterations (significantly fewer iterations than both other algoritms,
but also significantly lower restored PSNR).

\subsection{TV Denoising}

\begin{figure}
\centering
\begin{tabular}{c@{}c@{}c}
\epsfig{bb=90 90 410 410,clip=true,width=0.33\textwidth,figure=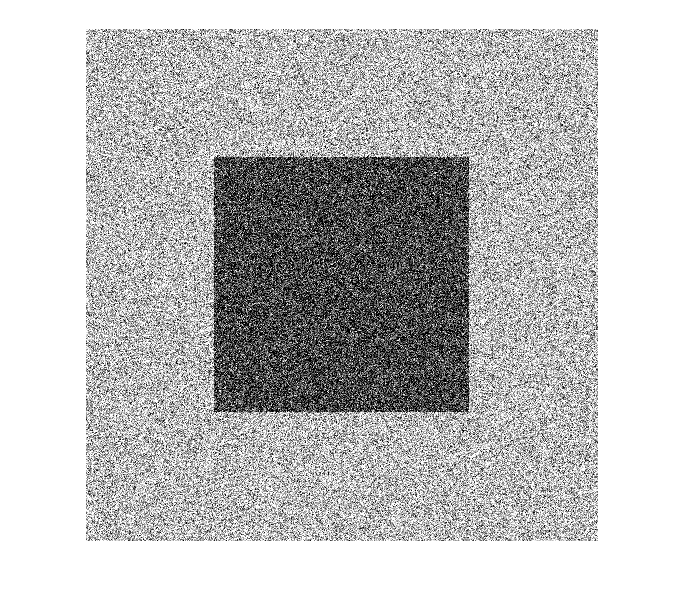} &
\epsfig{bb=90 90 410 410,clip=true,width=0.33\textwidth,figure=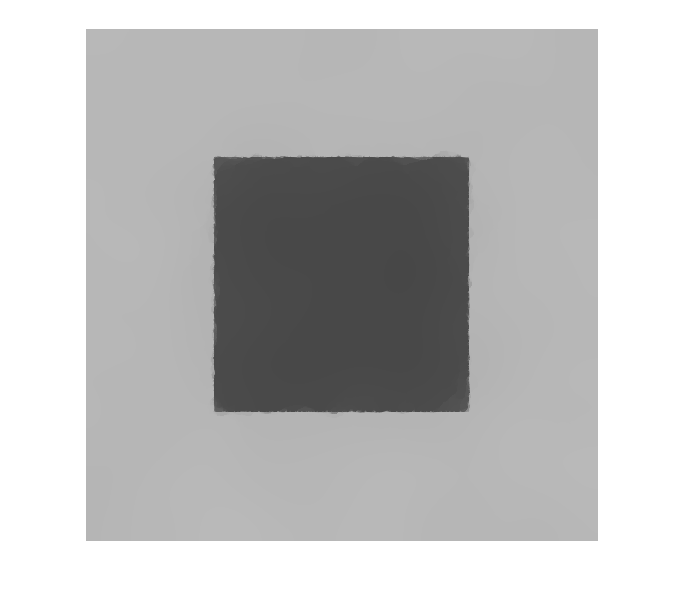} &
\epsfig{bb=90 90 410 410,clip=true,width=0.33\textwidth,figure=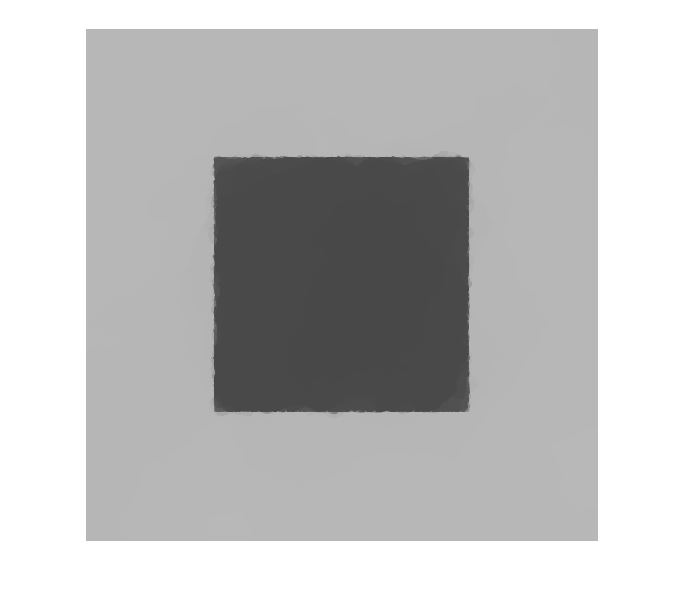} \\
Noisy Square & Split Bregman & PDE acceleration
\end{tabular}
\caption{Denoising of a synthetic image with total variation restoration with $\lambda=1000$ via (b) Split Bregman and (c) PDE acceleration. In PDE acceleration we used $\Delta t=\Delta x/2$ and $a=6\sqrt{\lambda}$.}
\label{fig:squareDN}
\end{figure}

\begin{figure}
\centering
\begin{tabular}{c@{}c}
\epsfig{trim=35 30 30 20,clip=true,width=0.5\textwidth,figure=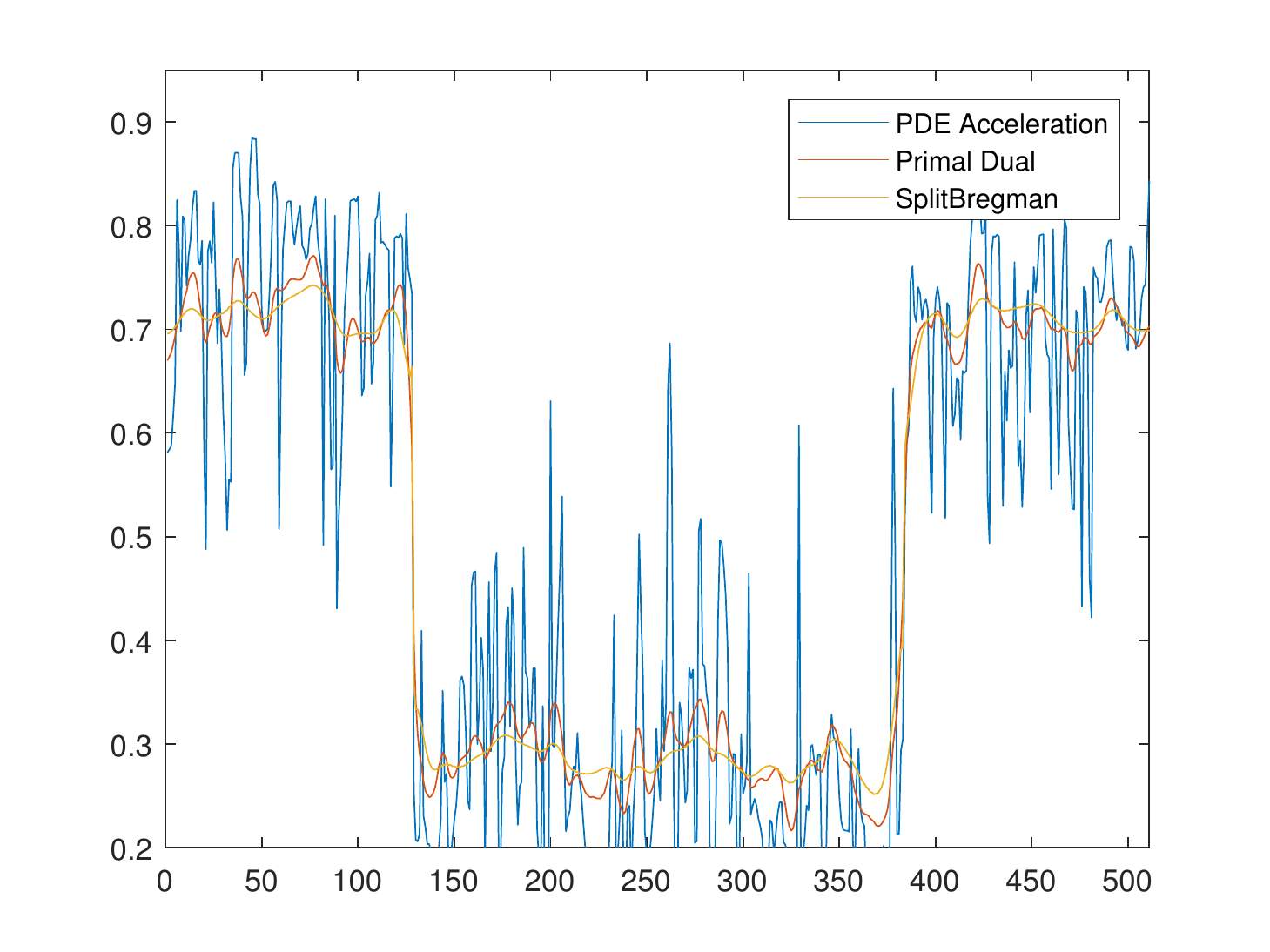} &
\epsfig{trim=35 30 30 20,clip=true,width=0.5\textwidth,figure=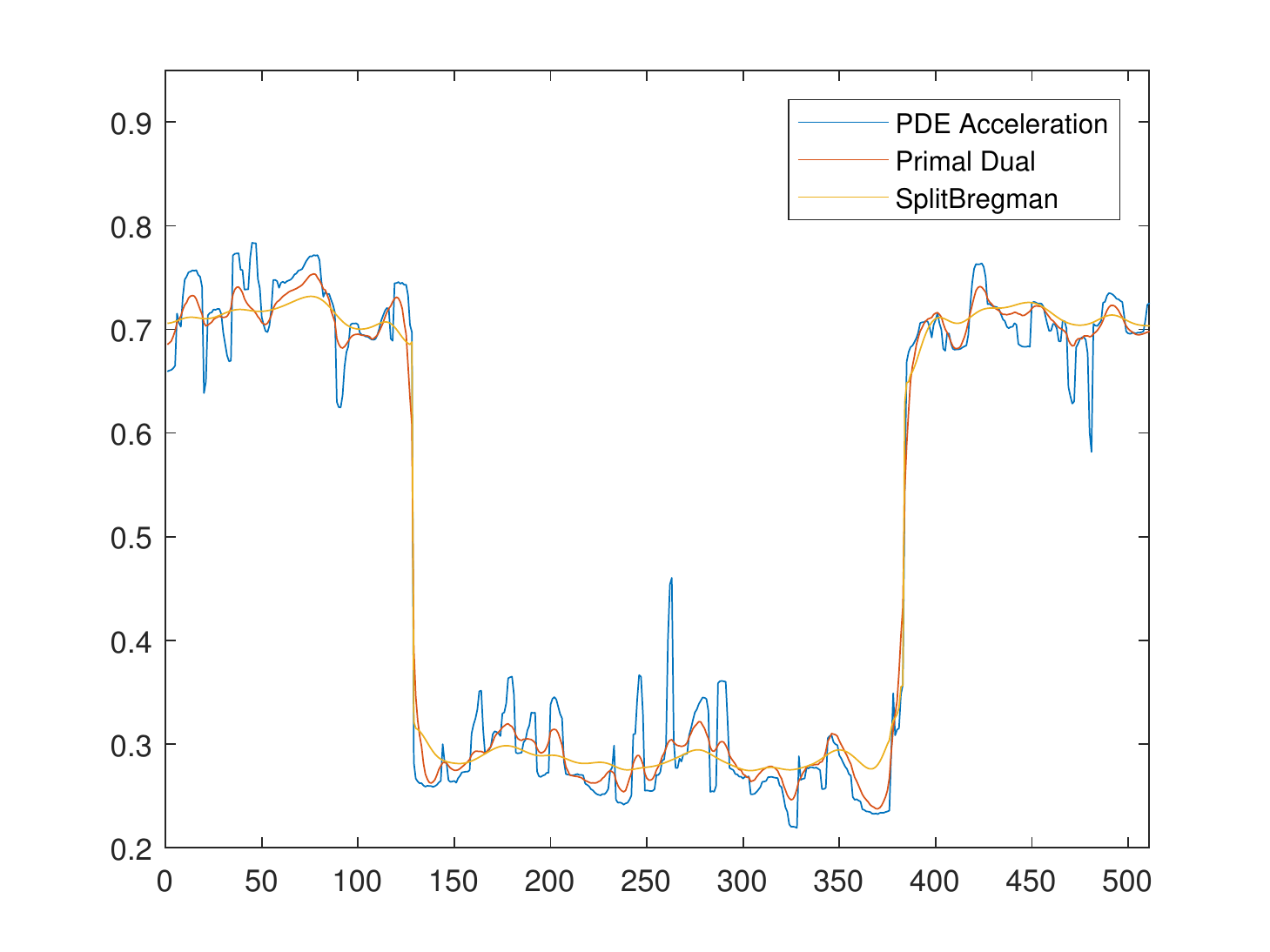} \vspace{-2mm}\\
$t=1$ & $t=2$ \\
\epsfig{trim=35 30 30 20,clip=true,width=0.5\textwidth,figure=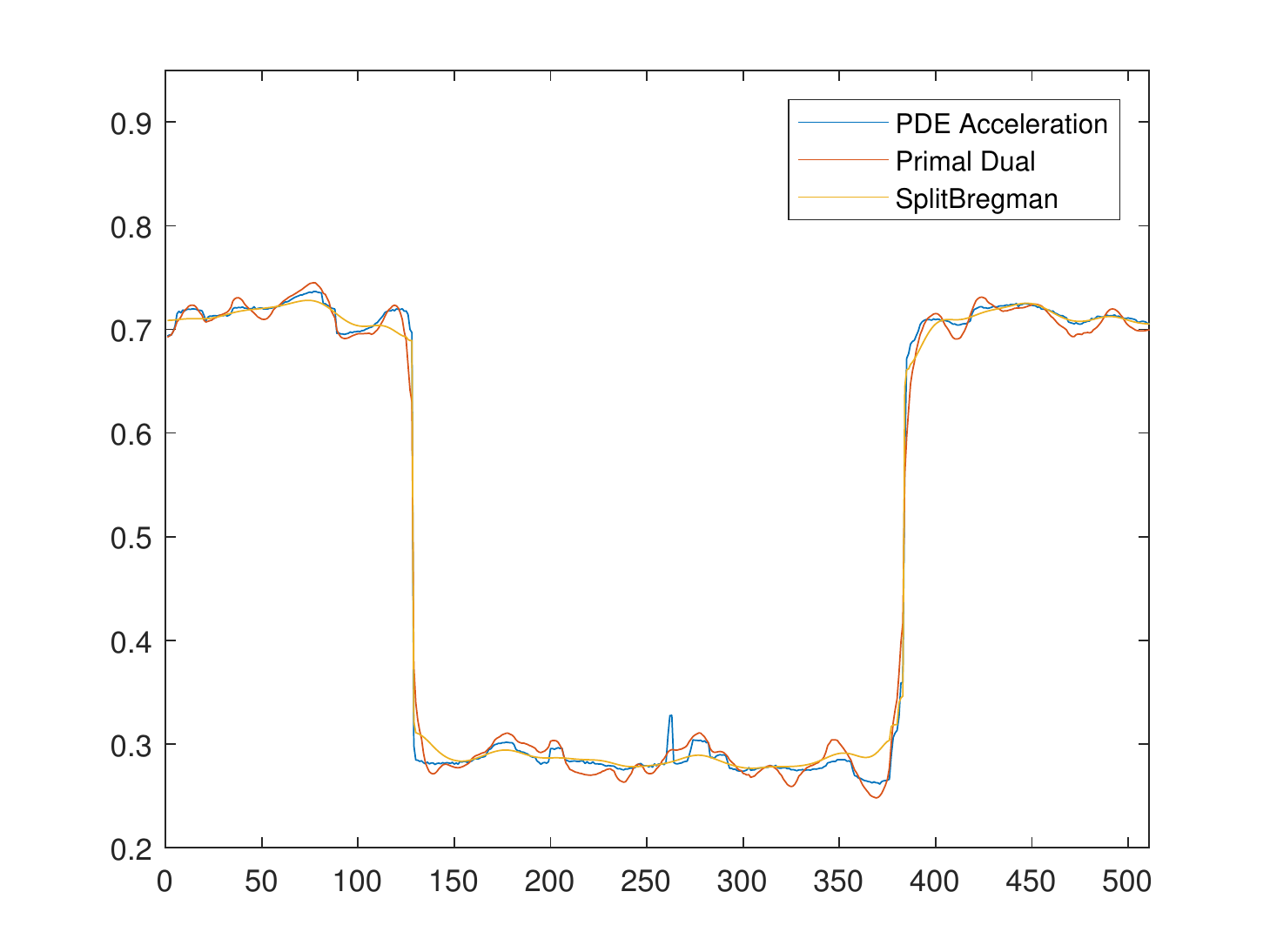} &
\epsfig{trim=35 30 30 20,clip=true,width=0.5\textwidth,figure=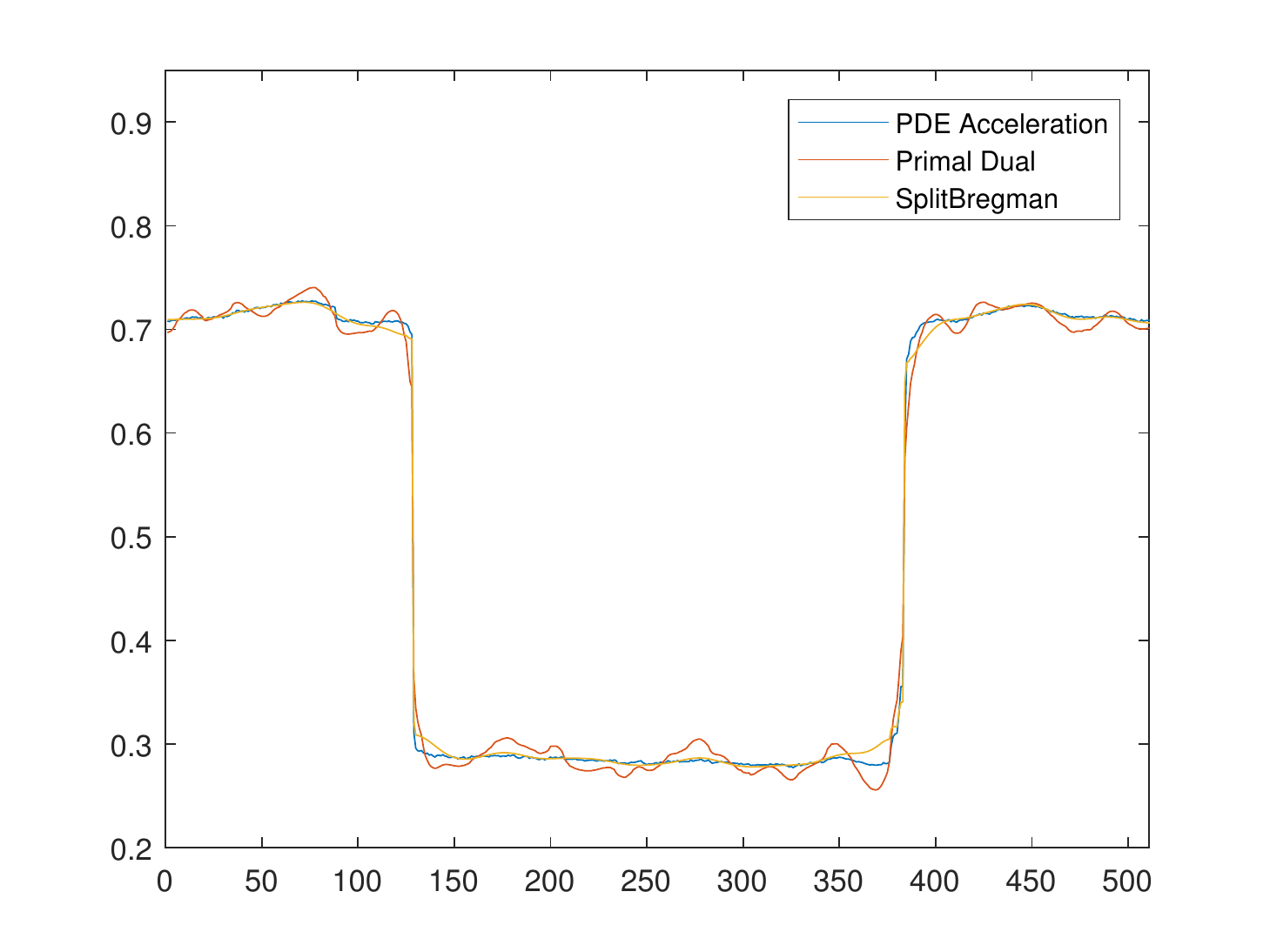} \vspace{-2mm}\\
$t=3$ & $t=4$
\end{tabular}
\caption{Comparison of PDE acceleration, Primal Dual, and Split Bregman algorithms for denoising a noisy square image. A one dimensional slice of the image is displayed at the same computation time for each algorithm. }
\label{fig:squareslice}
\end{figure}

\begin{figure}
\centering
\begin{tabular}{c@{}c}
\epsfig{trim=35 30 20 15,clip=true,width=0.50\textwidth,figure=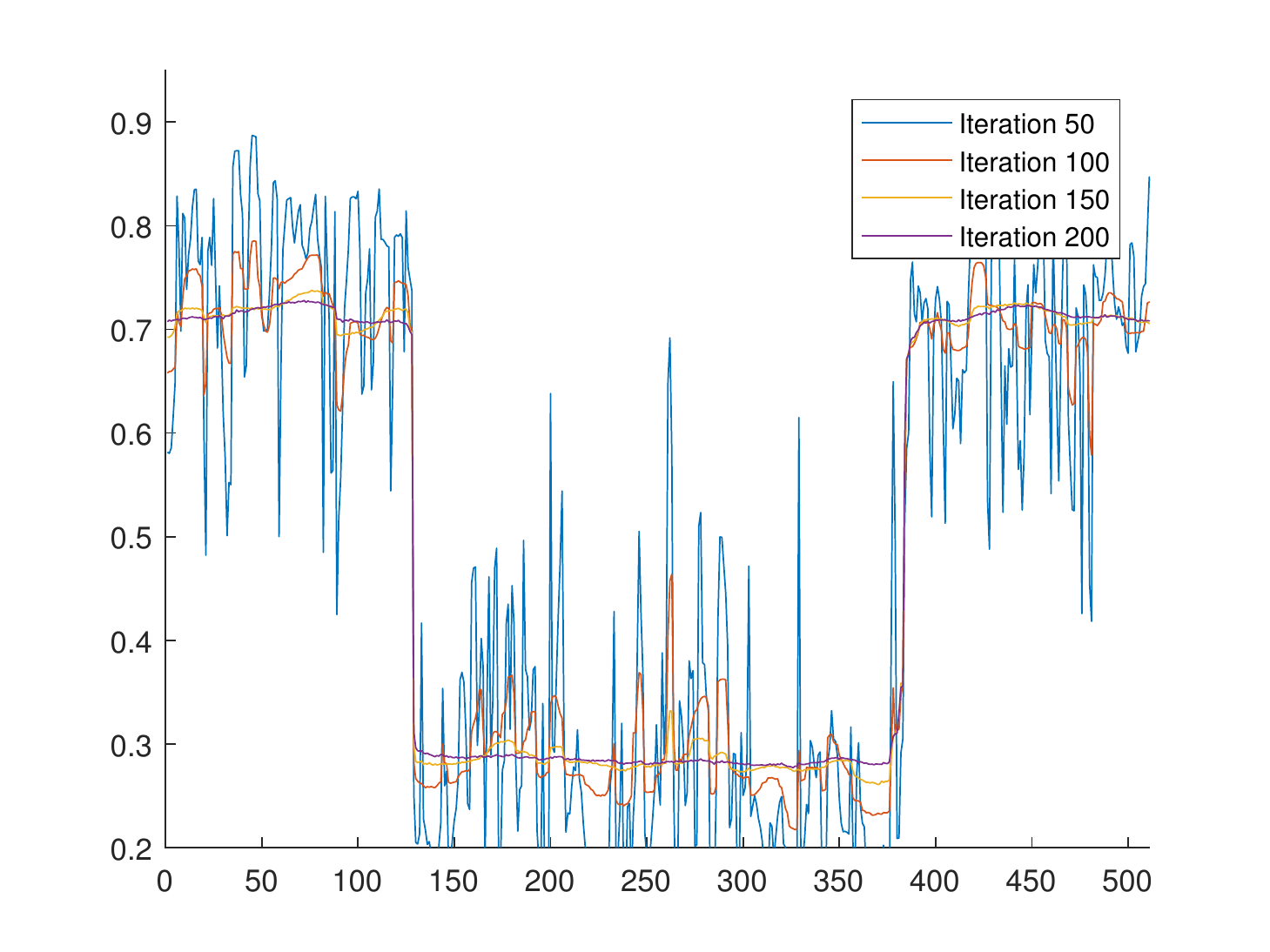} &
\epsfig{trim=35 30 20 15,clip=true,width=0.50\textwidth,figure=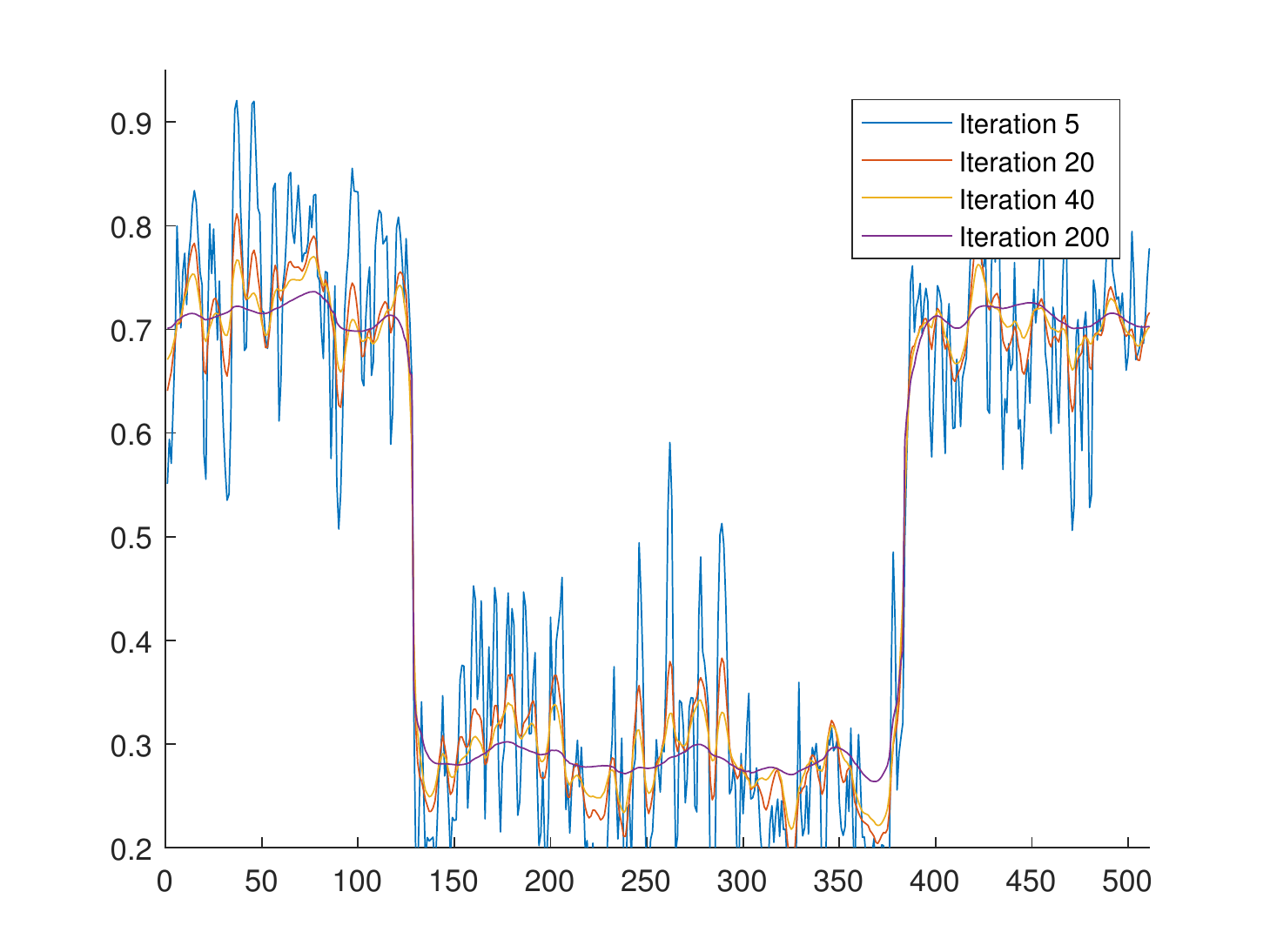} \\
PDE Acceleration & Primal Dual
\end{tabular}
\caption{Comparison of flows generated by (a) PDE Acceleration and (b) Primal Dual for solving the TV restoration problem on the noisy square image. Notice the edges are better preserved in PDE acceleration earlier in the flow. }
\label{fig:squaredemo}
\end{figure}

\begin{figure}
\centering
\begin{tabular}{c@{}c}
\epsfig{trim=35 30 20 15,clip=true,width=0.50\textwidth,figure=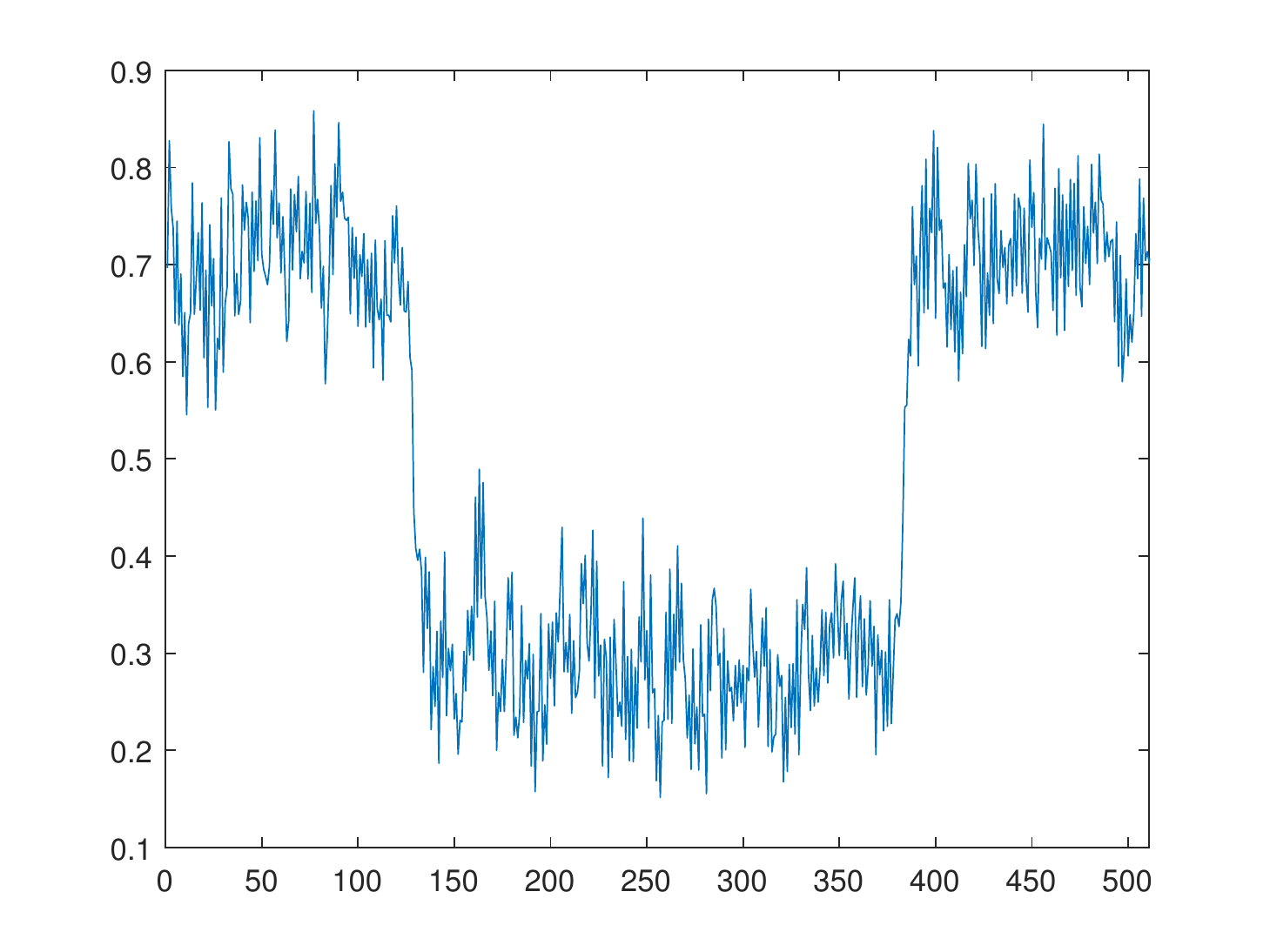} &
\epsfig{trim=35 30 20 15,clip=true,width=0.50\textwidth,figure=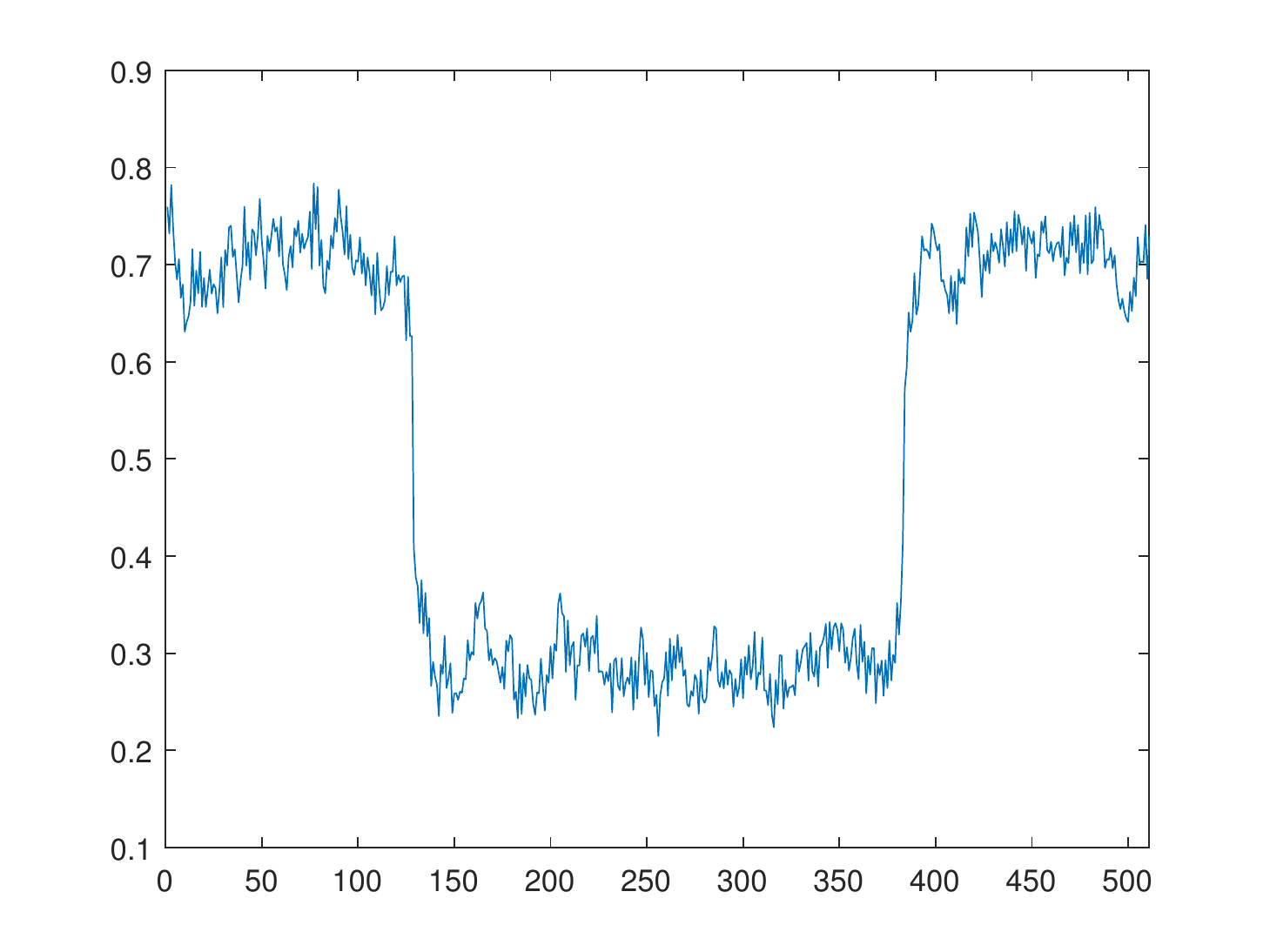} \\
$\Delta t=10\Delta x$ & $\Delta t=5\Delta x$ \\
\epsfig{trim=35 30 20 15,clip=true,width=0.50\textwidth,figure=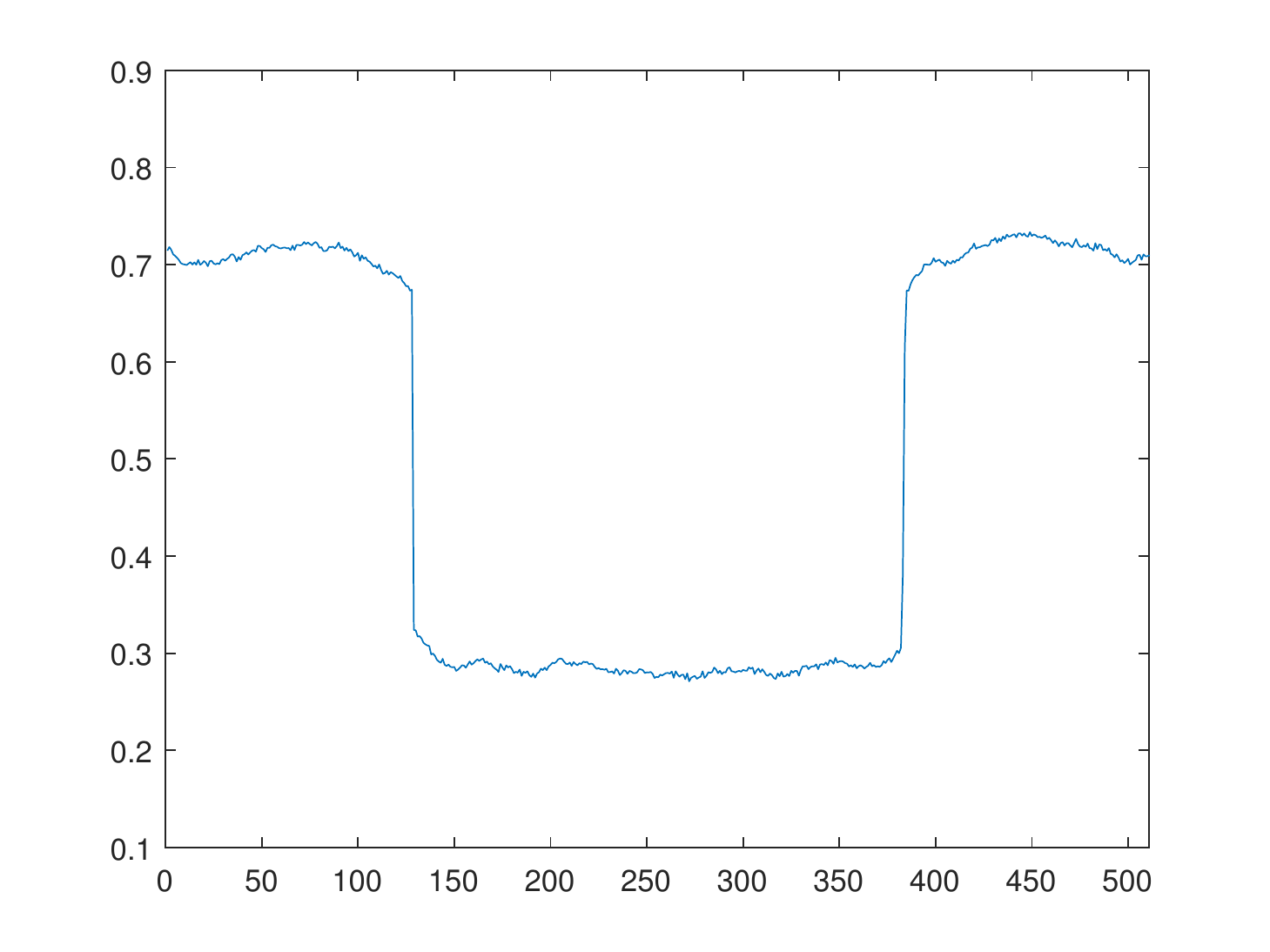} &
\epsfig{trim=35 30 20 15,clip=true,width=0.50\textwidth,figure=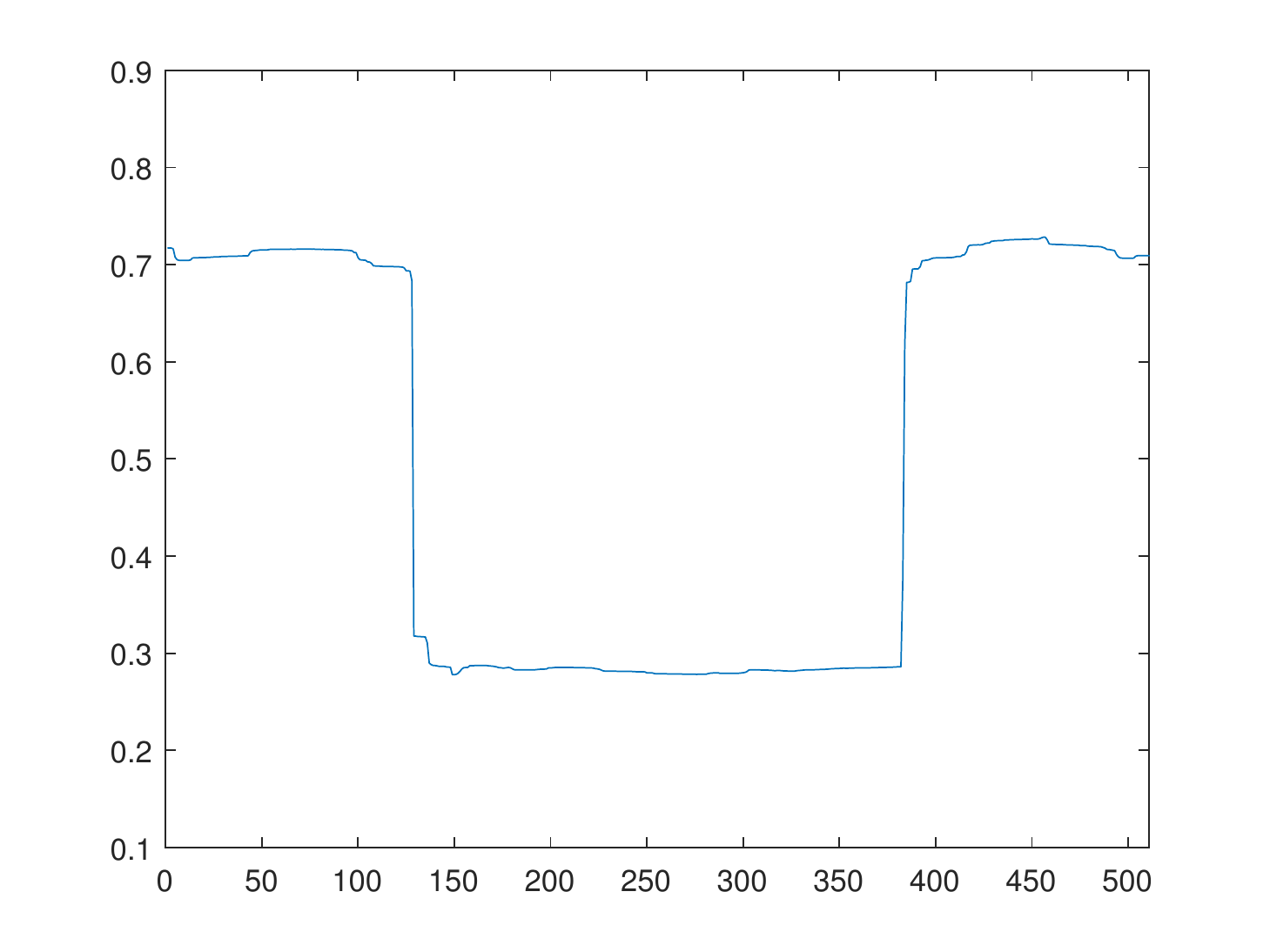} \\
$\Delta t=\Delta x$ & $\Delta t=\Delta x/10$
\end{tabular}
\caption{Comparison of steady state solutions for denoising the $2D$ noisy square for different time steps in PDE Acceleration. The scheme is stable in $L^\infty$ for a variety of time steps though we observe
$\Delta t\leq \Delta x$ is required to ensure the solution
is a reasonable denoising. }
\label{fig:cfl}
\end{figure}

\begin{figure}
\centering
\begin{tabular}{c@{}c}
\epsfig{trim=20 0 30 20,clip=true,width=0.50\textwidth,figure=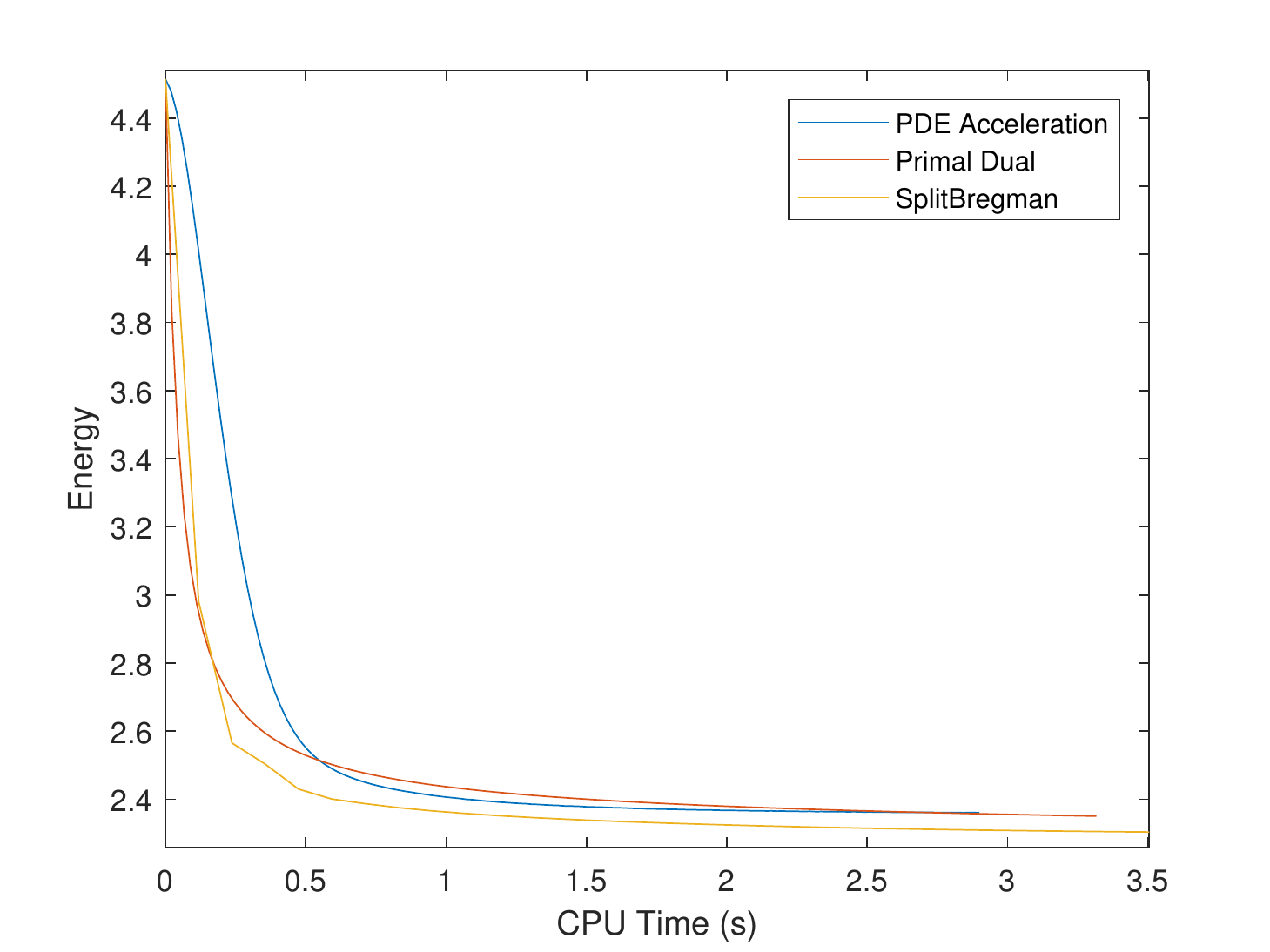} &
\epsfig{trim=20 0 30 20,clip=true,width=0.50\textwidth,figure=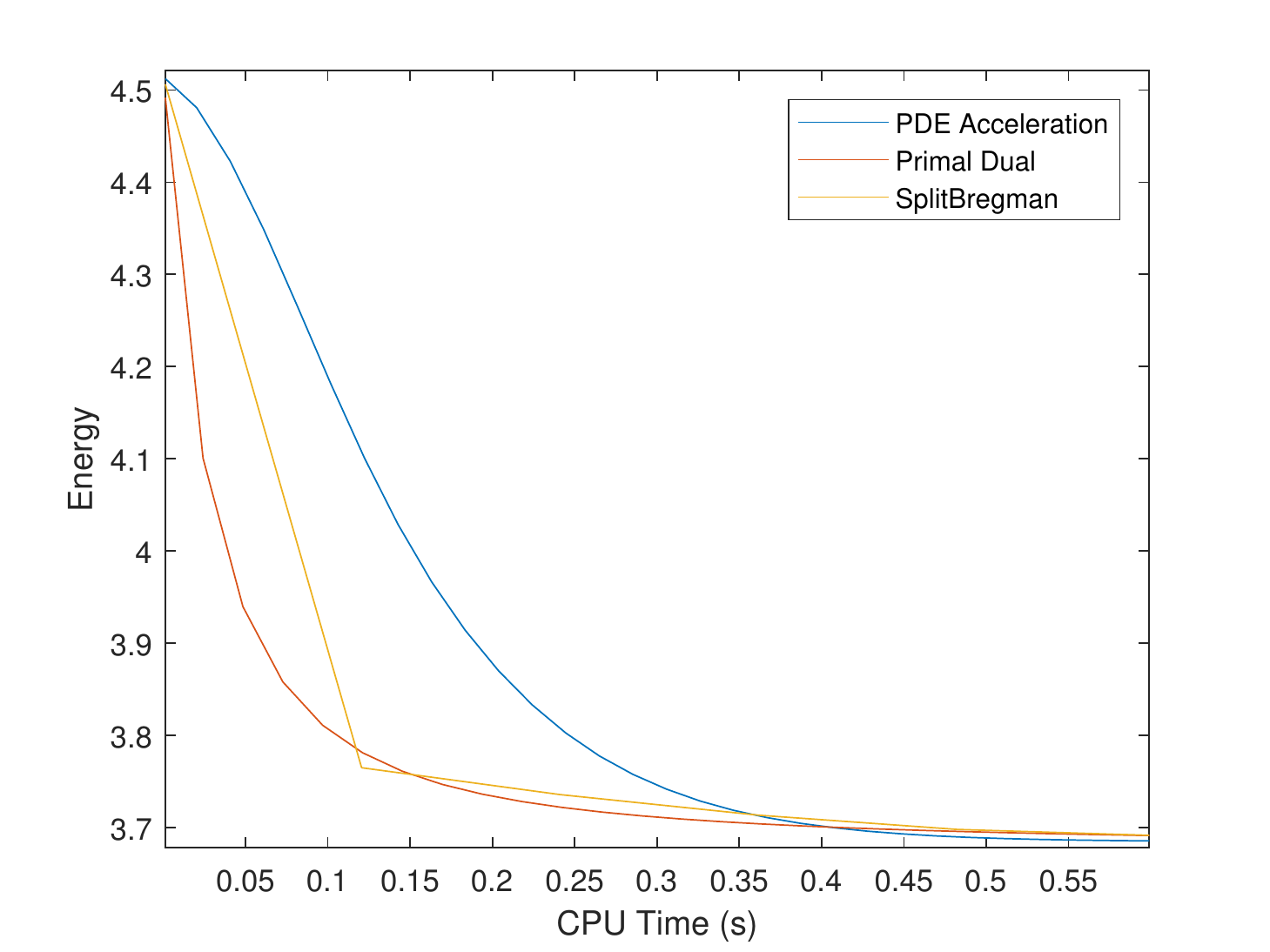} \\
$\lambda=1000$ & $\lambda=7000$
\end{tabular}
\caption{Comparison of logarithm of total energy versus CPU time for denoising the full $512\times 512$ Lenna image with PDE acceleration, Primal Dual, and Split Bregman. We used $a=6\sqrt{\lambda}$ for $\lambda=1000$ and $a=2\sqrt{\lambda}$ for $\lambda=7000$.}
\label{fig:E}
\end{figure}

\begin{figure}
\centering
\begin{tabular}{c@{}c}
\epsfig{trim=20 0 30 20,clip=true,width=0.50\textwidth,figure=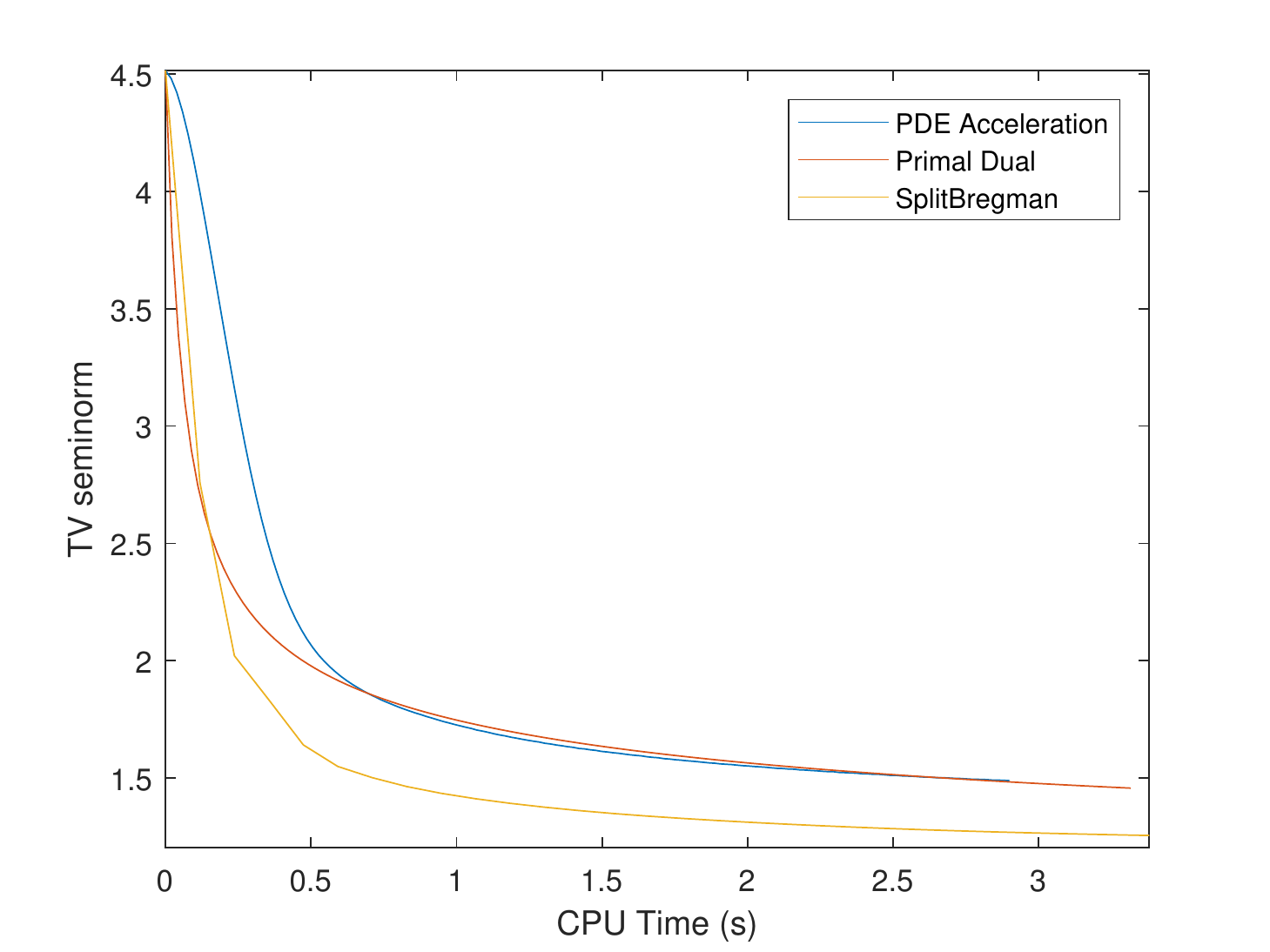} &
\epsfig{trim=20 0 30 20,clip=true,width=0.50\textwidth,figure=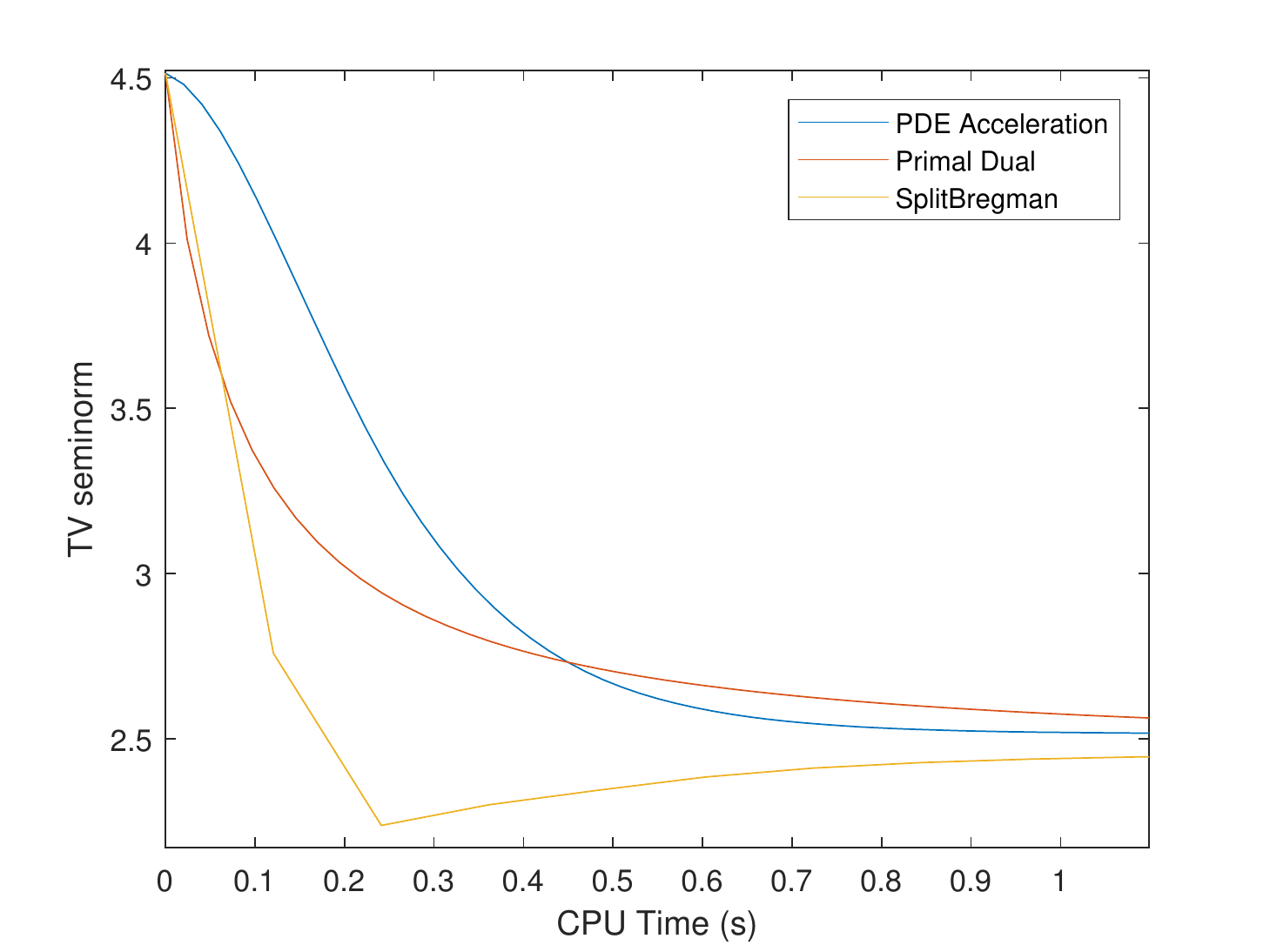} \\
$\lambda=1000$ & $\lambda=7000$
\end{tabular}
\caption{Comparison of logarithm of TV seminorm energy versus CPU time for denoising the full $512\times 512$ Lenna image with PDE acceleration, Primal Dual, and Split Bregman. }
\label{fig:E2}
\end{figure}

\begin{figure}
\centering
\begin{tabular}{c@{}c}
\epsfig{trim=20 0 30 20,clip=true,width=0.50\textwidth,figure=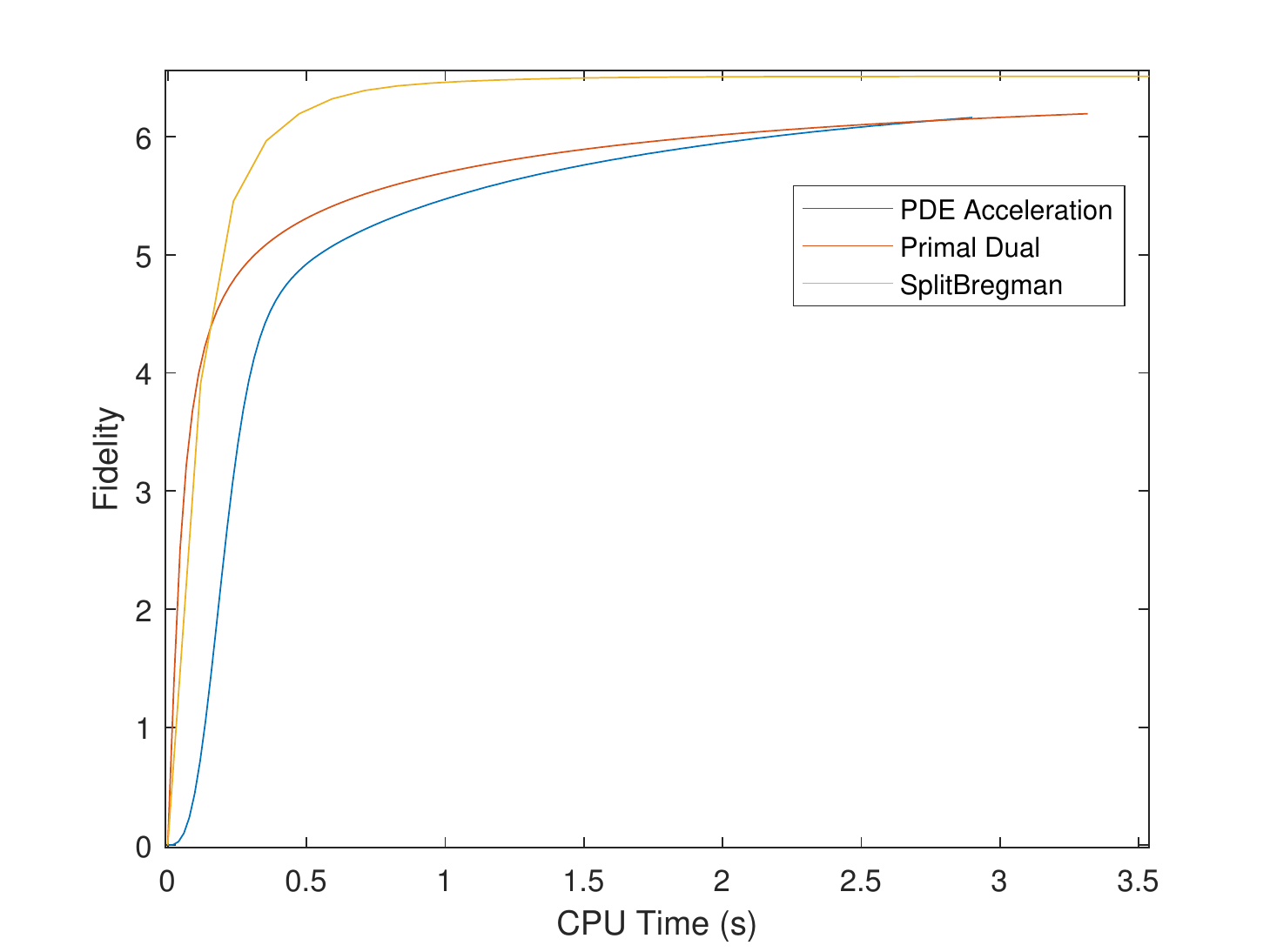} &
\epsfig{trim=20 0 30 20,clip=true,width=0.50\textwidth,figure=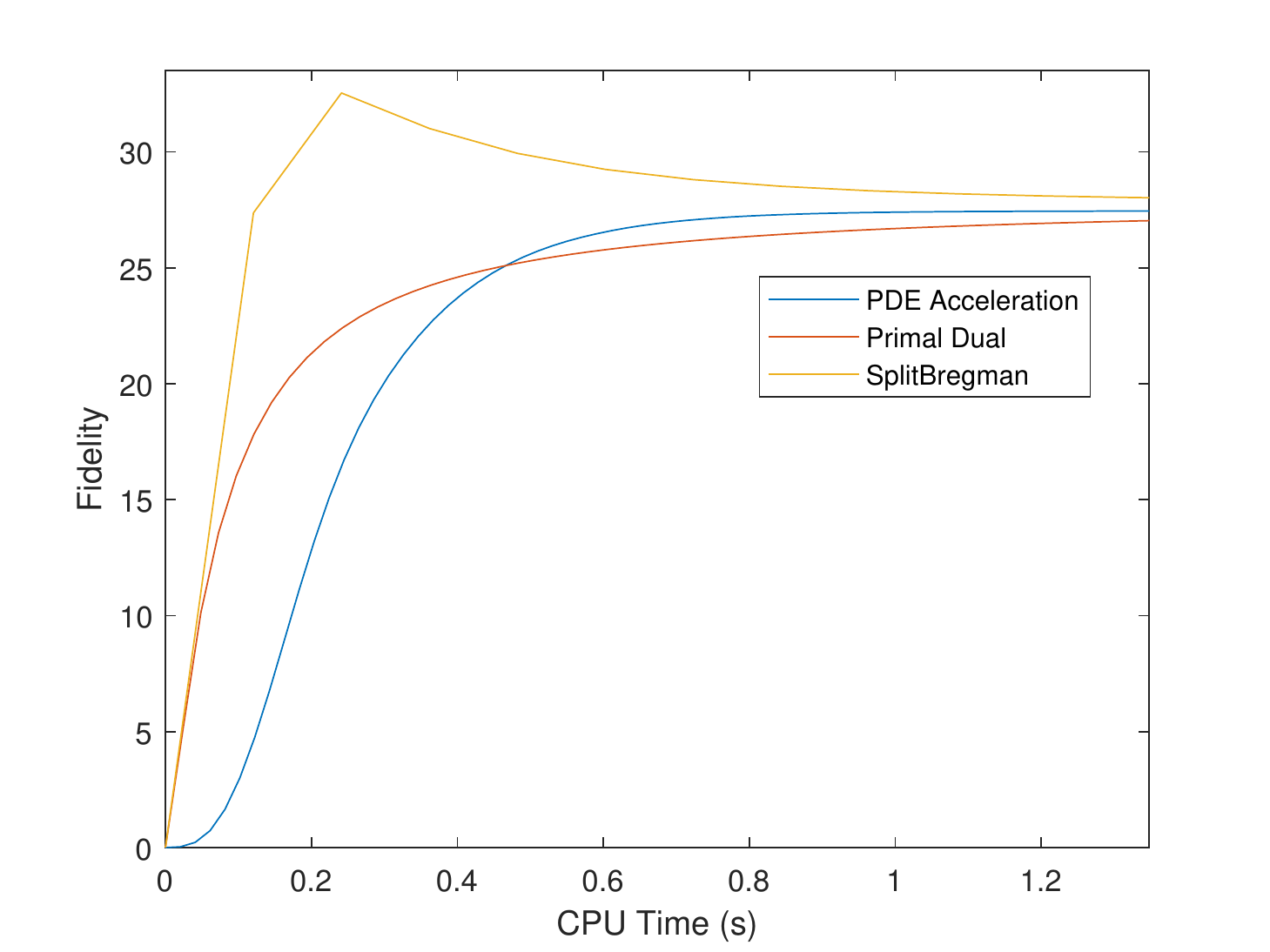} \\
$\lambda=1000$ & $\lambda=7000$
\end{tabular}
\caption{Comparison of fidelity energy versus CPU time for denoising the full $512\times 512$ Lenna image with PDE acceleration, Primal Dual, and Split Bregman. }
\label{fig:E3}
\end{figure}

\begin{figure}
\centering
\begin{tabular}{c@{}c@{}c@{}c@{}c}
\epsfig{trim=120 90 300 50,clip=true,width=0.20\textwidth,figure=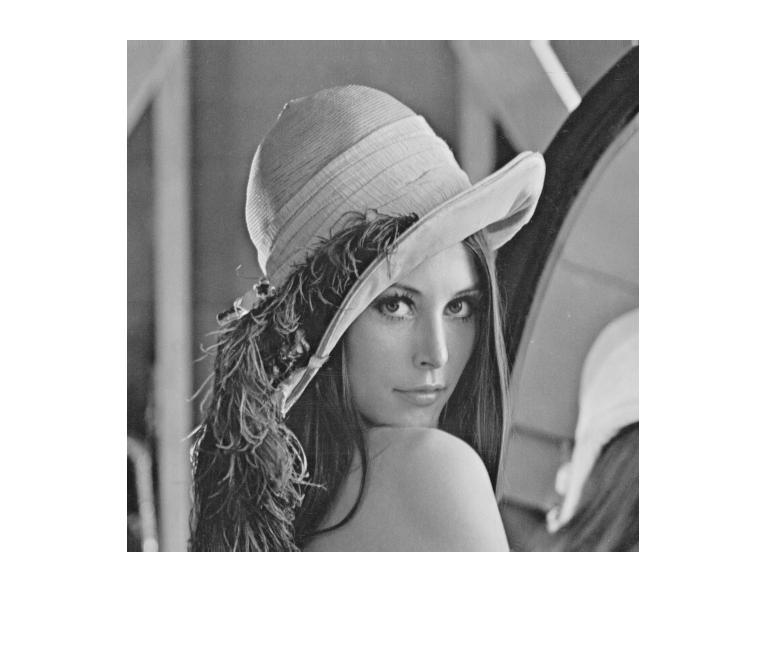} &
\epsfig{trim=120 90 300 50,clip=true,width=0.20\textwidth,figure=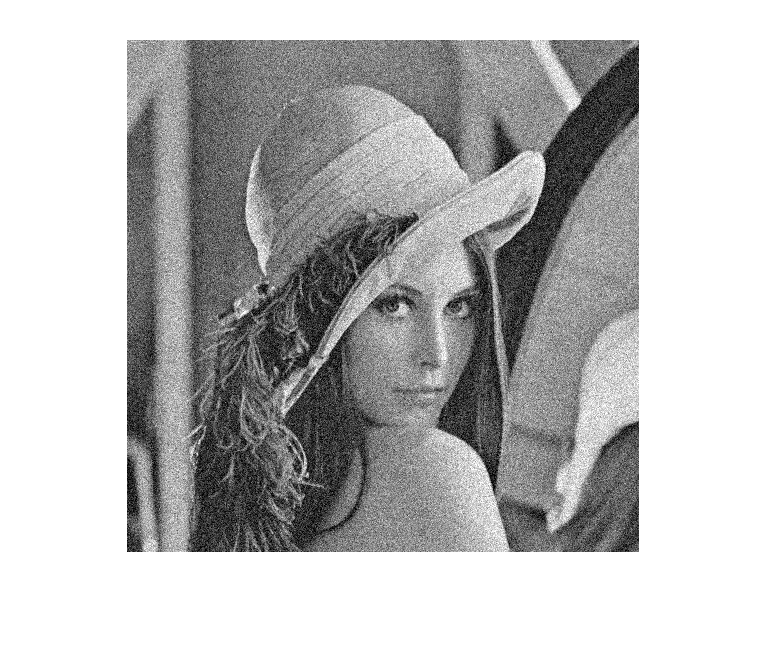} &
\epsfig{trim=120 90 300 50,clip=true,width=0.20\textwidth,figure=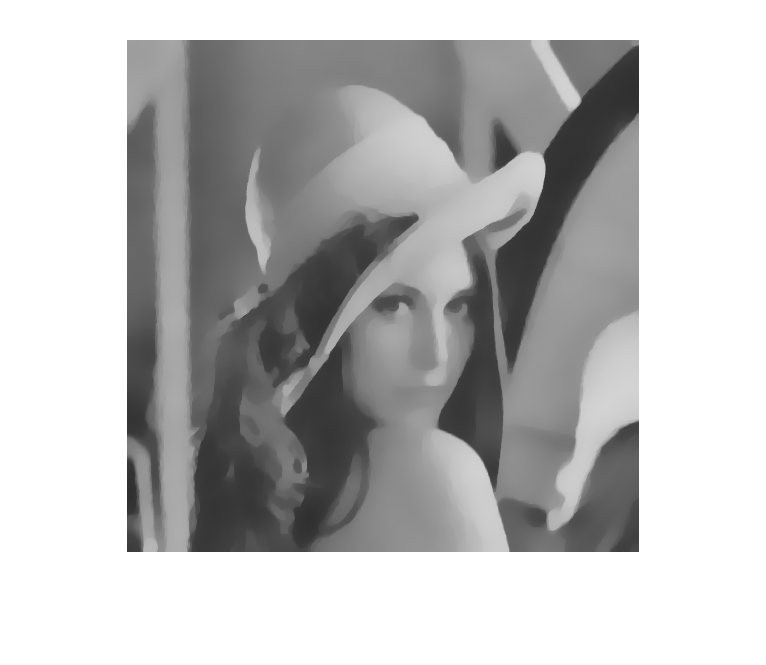} &
\epsfig{trim=120 90 300 50,clip=true,width=0.20\textwidth,figure=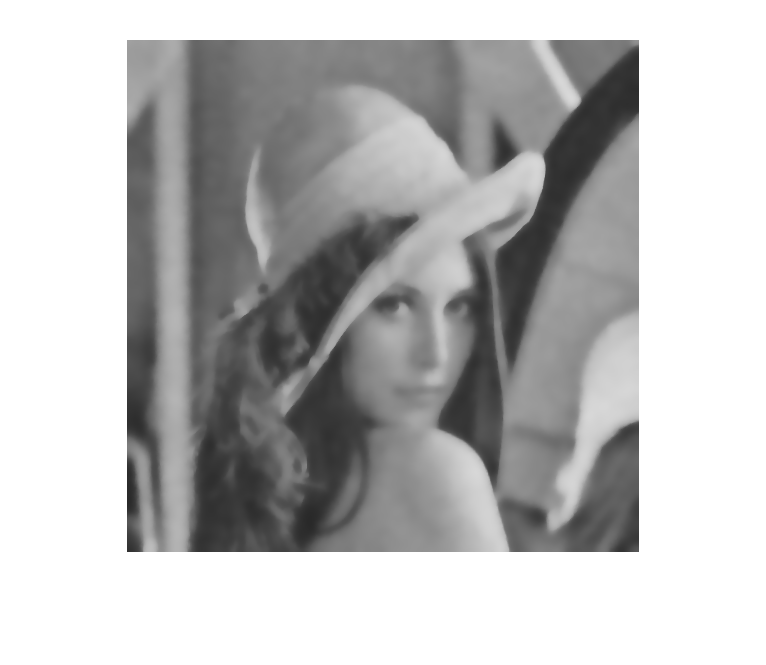} &
\epsfig{trim=120 90 300 50,clip=true,width=0.20\textwidth,figure=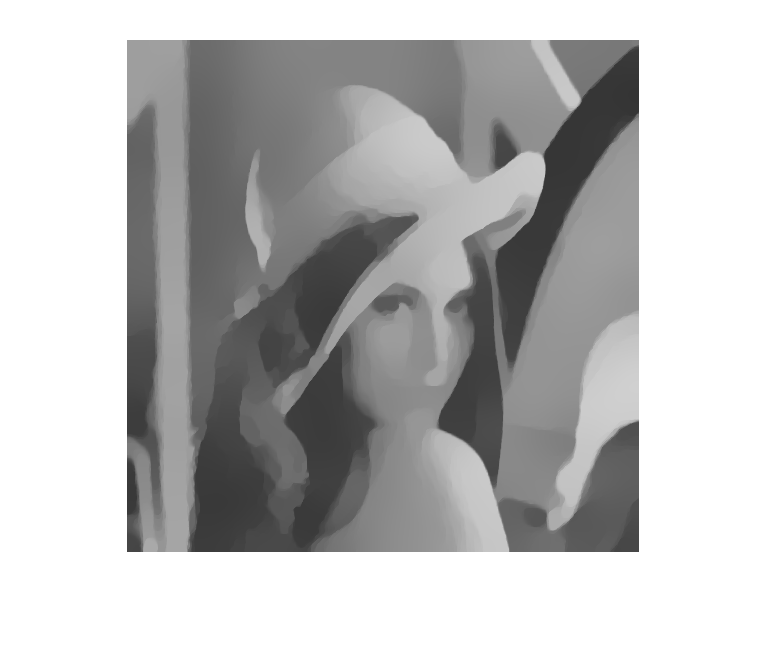} \\
Original & Noisy & PDE acceleration & Primal Dual & Split Bregman
\end{tabular}
\caption{Comparison of PDE acceleration, Primal Dual, and Split Bregman for TV restoration of a noisy Lenna image with $\lambda=1000$. Each algorithm was run for 150 iterations, which took 2.7 seconds for PDE acceleration, 3.3 seconds for Primal Dual, and 28 seconds for Split Bregman.}
\label{fig:lenna1}
\captionsetup[subfigure]{labelformat=parens}
\end{figure}

\begin{figure}
\captionsetup[subfigure]{labelformat=empty}
\centering
\begin{tabular}{c@{}c@{}c@{}c@{}c}
\epsfig{trim=120 90 300 50,clip=true,width=0.20\textwidth,figure=lenna} &
\epsfig{trim=120 90 300 50,clip=true,width=0.20\textwidth,figure=lenna_noisy} &
\epsfig{trim=120 90 300 50,clip=true,width=0.20\textwidth,figure=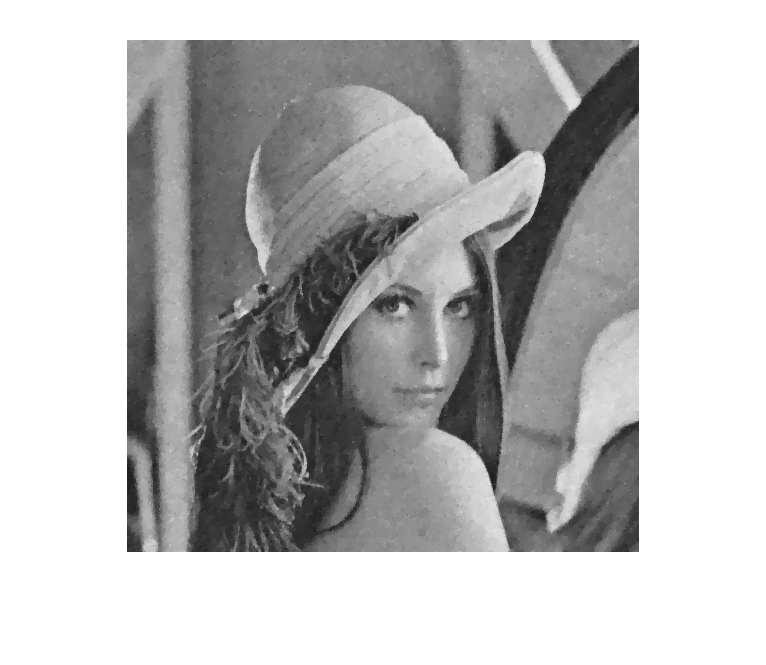} &
\epsfig{trim=120 90 300 50,clip=true,width=0.20\textwidth,figure=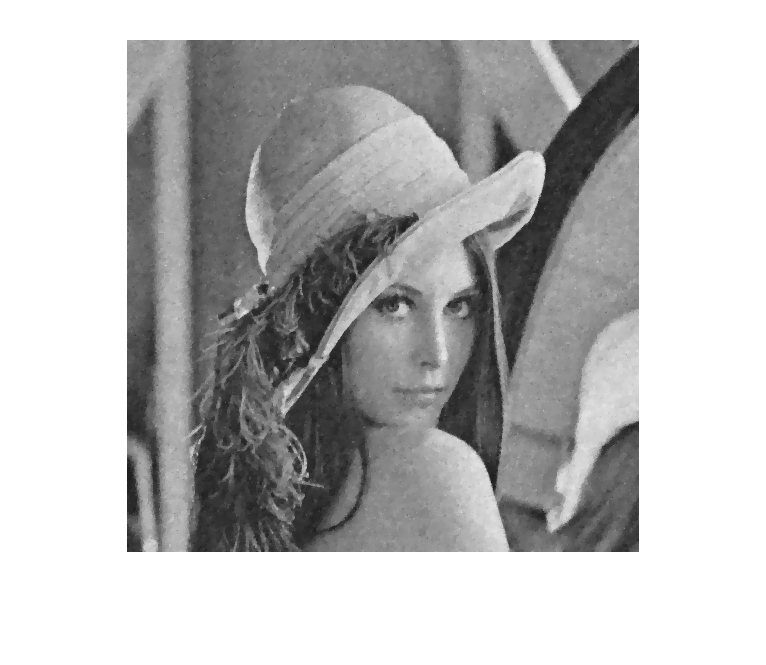} &
\epsfig{trim=120 90 300 50,clip=true,width=0.20\textwidth,figure=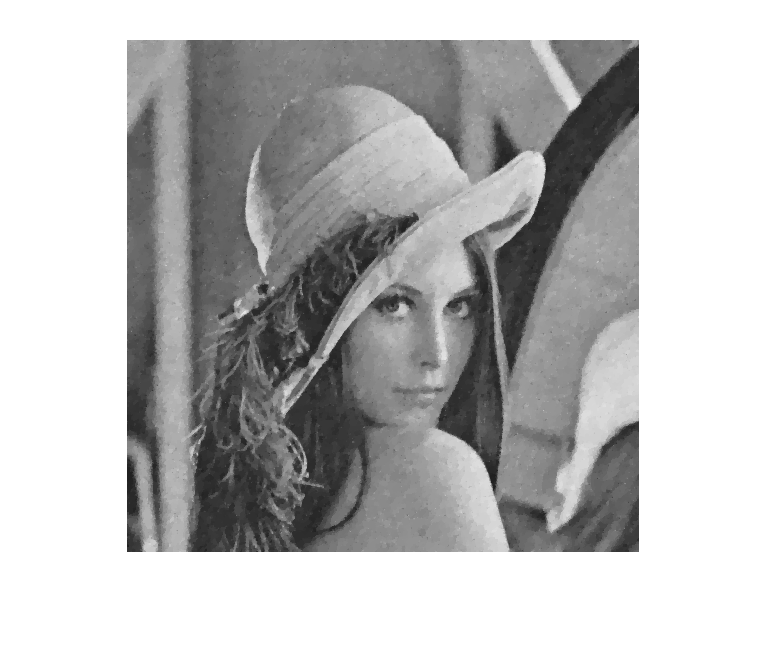} \\
Original & Noisy & PDE acceleration & Primal Dual & Split Bregman
\end{tabular}
\caption{Comparison of PDE acceleration, Primal Dual, and Split Bregman for TV restoration of a noisy Lenna image with $\lambda=7000$. Each algorithm was run for 50 iterations, which took 0.85 seconds for PDE acceleration, 1.12 seconds for Primal Dual, and 10.4 seconds for Split Bregman.}
\label{fig:lenna2}
\captionsetup[subfigure]{labelformat=parens}
\end{figure}

We now consider the problem of Total Variation (TV) restoration, which
has a long history in image processing \cite{rudin1992nonlinear}.
The TV denoising problem
corresponds to minimizing \eqref{eq:TV-energy} in the absence
of a kernel $K$ via the accelerated PDE \eqref{eq:TV-PDE}. In this case
state of the art approaches include primal dual methods
\cite{chambolle2011first} and the split Bregman method
\cite{goldstein2009split}.

We again use the first order explicit scheme \eqref{eq:TV-forward},
while discretizing the spatial gradient and divergence separately (using forward
differences for the gradient and backward differences for the divergence),
and homogeneous Neumann boundary conditions.
Numerically, we set $\nabla u/|\nabla u|=0$
whenever $\nabla u=0$, so no regularization is required, though we
rarely encounter numerical gradients that are identically zero. This
choice of discretization makes the discrete divergence the exact numerical
adjoint of the discrete gradient.

We first consider a noisy square image, with dark region $u=0.25$
and light region $u=0.75$ with additive Gaussian noise with standard
deviation $\sigma=0.3$. Figure \ref{fig:squareDN} shows the noisy
square and the total variation denoising with the Split Bregman algorithm
and PDE acceleration. We compare PDE acceleration, Primal Dual, and
Split Bregman on slices of the image at similar computation times
in Figure \ref{fig:squareslice}. Notice the Primal Dual algorithm
blurs the edges slightly at first, and they are restored only late
in the flow (at $t=4$ Primal dual has not yet converged). The PDE
acceleration algorithm does a better job preserving edges (they are
never blurred) compared to Primal Dual, and is slightly better than
Split Bregman at preserving edges by time $t=4$.

In the example above we took $\Delta t=\Delta x/2$ for simplicity.
Corroborating our analysis in Section \ref{sec:TV}, this explicit numerical
scheme \eqref{eq:TV-forward} behaves stably
in $L^{\infty}$ in our experiments, meaning the solutions remain
bounded in $L^{\infty}$ for all time, even for larger time steps which
still satisfy the necessary conditions \eqref{eq:TV-gradient-CFL},
\eqref{eq:TV-forward-CFL}, \eqref{eq:TV-central-CFL} or \eqref{eq:TV-semi-CFL}.
For such larger time steps, though, we find the flow does not fully converge,
yet remains stable via the nonlinear effect discussed in Section \ref{sec:TV},
but instead tends to an oscillatory steady state. Figure \ref{fig:cfl}
shows a snapshot of the steady state for various values of the time
step $\Delta t$. For $\Delta t\leq \Delta x$ the steady state is a
reasonable denoising,
hence we choose $\Delta t=\Delta x$ or $\Delta t=\Delta x/2$
in most of this paper. Note that this closely matches the suggested time
step bound in \eqref{eq:TV-central} for a quantization level of 1/255,
given the other parameters utilized here, which would come out to
$\Delta t\le1.189\Delta x$.

Figure \ref{fig:E} compares the energy decay against CPU time for
denoising the Lenna image with PDE acceleration, Primal Dual, and
Split Bregman algorithms. The noise is additive zero mean Gaussian
noise with standard deviation $\sigma=0.1$ and the images take values
in the interval $[0,1]$. We note in Figures \ref{fig:lenna1} and
\ref{fig:lenna2} that PDE acceleration appears to yield a better
quality image for the same energy level compared to primal dual.

\section{Conclusion} \label{sec:conc}

We employed the novel framework of \emph{PDE Acceleration},
based on momentum  methods such as Nesterov and Polyak's heavy ball method,
to calculus of variations problems defined for general
functions on $\R^n$. The result was a very general set of
\emph{accelerated PDE's}
whose simple discretizations efficiently solve the
the resulting class of optimization problems. We further analyzed their
use in regularized inversion problems, where gradient descent diffusion
equations get replaced by nonlinear wave equations within the framework
of PDE acceleration, with far more generous discrete time step conditions.

We presented results of experiments on image processing
problems including Beltrami regularized denoising and inpainting,
and total variation (TV) regularized denoising and deblurring. In
all cases, we can achieve state of the art results with very simple
algorithms; indeed, the PDE acceleration update is a simple explicit
forward Euler update of a nonlinear wave equation. Future work will
focus on problems such as TV inpainting, where there is no fidelity,
how to choose the damping parameter adaptively to further accelerate
convergence, and applications to other problems in computer vision,
such as Chan-Vese active contours \cite{chan2001active}.


\end{document}